\documentclass[11pt]{article}
\usepackage{lmodern}
\usepackage{fullpage}
\usepackage{authblk}
\usepackage{mathrsfs}
\usepackage{amssymb}
\usepackage{amsthm}
\usepackage{cjhebrew}
\usepackage{bbm}
\usepackage{bbold}
\usepackage{stmaryrd}
\usepackage{setspace}
\usepackage{xspace}
\usepackage{enumerate}
\usepackage{enumitem}
\usepackage{mathtools}
\usepackage{thm-restate}
\usepackage[hidelinks]{hyperref}
\usepackage{cleveref}
\usepackage{xcolor}

\usepackage{tocloft}
\emergencystretch=\maxdimen
\allowdisplaybreaks

\declaretheorem[sibling=basecase]{theorem, proposition, lemma, corollary, claim, conjecture, observation}
\declaretheorem[style=definition, sibling=basecase]{definition, convention, notation, example, remark, remarks, question, assumption, comment}
\declaretheorem[name=Theorem, sibling=basecase, numbered=no]{theorem*}

\newlist{thmlist}{enumerate}{1}  
\setlist[thmlist]{label=\roman{thmlisti}., ref=\thetheorem.(\roman{thmlisti}),noitemsep} 

\addtotheorempostheadhook[theorem]{\crefalias{thmlisti}{theorem}}
\addtotheorempostheadhook[lemma]{\crefalias{thmlisti}{lemma}}
\addtotheorempostheadhook[proposition]{\crefalias{thmlisti}{proposition}}

\Crefname{theorem}{Theorem}{Theorems}
\Crefname{lemma}{Lemma}{Lemmas}
\Crefname{proposition}{Proposition}{Propositions}

\newcommand{\tf}[1]{\textnormal{\textsc{#1}}\xspace}
\newcommand{\ff}[1]{\textnormal{\textrm{#1}}\xspace}

\DeclareMathOperator{\ZFepsilon}{\tf{ZF}_{\varepsilon}}
\DeclareMathOperator{\GBepsilon}{\tf{GB}_{\varepsilon}}
\newcommand{\ZFepsilonL}[1]{\tf{ZF}_{\varepsilon}^{#1}}

\newcommand\cnot[1]{\mathrel{\ooalign{\hfil$#1$\hfil\cr\hfil$/$\hfil\cr}}}
\newcommand{\rlzin}{\mathop{\varepsilon}}
\newcommand{\notrlzin}{\cnot{\varepsilon}}
\newcommand{\falsity}[1]{\lVert #1 \rVert}
\newcommand{\verity}[1]{\lvert #1 \rvert}
\newcommand{\Perp}{\mathbin{\text{$\bot\mkern-11mu\bot$}}}
\newcommand{\point}{\boldsymbol{.}}

\newcommand{\stackapp}{\point}
\newcommand{\lambdaapp}{\point}

\DeclareMathOperator{\divline}{\hspace{0.05cm}\mid\hspace{0.05cm}}

\newcommand{\rlzfont}[1]{\textnormal{\textbf{#1}}}
\newcommand{\saverlz}[1]{\ff{k}_{#1}}
\newcommand{\rlzstr}{\tf{N}}
\newcommand{\rlzmodel}{\mathcal{N}}
\newcommand{\rlzset}{\mathcal{R}}
\newcommand{\fullname}[1]{\text{\Large \cjRL{r}}({#1})} 
\newcommand{\cjgimel}{\text{\Large \cjRL{g}}}
\newcommand{\identity}{\rlzfont{I}}
\newcommand{\cc}{\textnormal{\ff{cc}}}

\newcommand{\pointwise}{\textnormal{``}}

\newcommand{\lquote}{\text{``}}
\newcommand{\rquote}{\text{''}}

\newcommand{\sing}{\textnormal{\texttt{sng}}}
\newcommand{\up}{\textnormal{\texttt{up}}}
\newcommand{\op}{\textnormal{\texttt{op}}}
\newcommand{\PlusOne}{\widetilde{\ff{succ}}}

\newcommand{\app}[2]{#1 #2}

\newcommand{\lapp}[2]{(#1)#2}
\newcommand{\rapp}[2]{#1(#2)}
\newcommand{\twoapp}[2]{(#1)(#2)}

\newcommand{\biglapp}[2]{\big(#1\big)#2}
\newcommand{\bigrapp}[2]{#1\big(#2\big)}
\newcommand{\bigtwoapp}[2]{\big(#1\big)\big(#2\big)}

\newcommand{\Bigtwoapp}[2]{\Big(#1\Big)\Big(#2\Big)}

\newcommand{\biggrapp}[2]{#1\bigg(#2\bigg)}

\newcommand{\Bigglapp}[2]{\Bigg(#1\Bigg)#2}

\newcommand{\lift}[1]{\tilde{#1}}
\newcommand{\imp}{\rightarrow}

\newcommand{\rlzXSubsetX}{\rlzfont{w}_0}
\newcommand{\rlzXSimeqX}{\rlzfont{w}_1}
\newcommand{\rlzXNotInX}{\rlzfont{w}_2}
\newcommand{\rlzSubsetTrans}{\rlzfont{w}_3}
\newcommand{\rlzInIsNotSimeq}{\rlzfont{w}_4}
\newcommand{\rlzNotNot}{\rlzfont{w}_5}
\newcommand{\rlzXNotSimeqY}{\rlzfont{w}_6}

\newcommand{\lamtermOne}{u}
\newcommand{\lamtermTwo}{v}
\newcommand{\lamtermThree}{w}
\newcommand{\lamAbst}[1]{\lambda #1 \lambdaapp}
\newcommand{\lamAbstOne}{\lamAbst{\lamtermOne}}
\newcommand{\lamAbstTwo}{\lamAbst{\lamtermTwo}}
\newcommand{\lamAbstThree}{\lamAbst{\lamtermThree}}

\title{Realizing the totally unordered structure of ordinals}
\date{}
\author{Laura Fontanella and Richard Matthews}
\affil{Univ. Paris Est Créteil, LACL, F-94010}

\newcommand{\Rcomment}[1]{{\color{black}{#1}}}
\newcommand{\Lcomment}[1]{{\color{black}{#1}}}

\begin{document}

\maketitle

\begin{abstract}
    We present tools for analysing ordinals in realizability models of \Rcomment{classical} set theory \Lcomment{built using Krivine's technique for realizability.} \Lcomment{This method uses a conservative extension of ZF known as $\ZFepsilon,$ where two membership relations co-exist, the usual one denoted $\in$ and a stricter one denoted $\rlzin$ that does not satisfy the axiom of extensionality; accordingly we have two equality relations, the extensional one $\simeq$ and the strict identity $=$ referring to sets that satisfy the same formulas.} We define recursive names using an operator that we call \emph{reish} and denote $\text{\Large \cjRL{r}}$, and we show that the class of recursive \Rcomment{names} for ordinals \Rcomment{coincides} extensionally with the class of ordinals of realizability models. We show that $\fullname \omega$ is extensionally \Rcomment{equal} to omega in any realizability model, thus recursive names provide a useful tool for computing $\omega$ in realizability models. We show that on the contrary $\rlzin$-totally ordered sets 
    do not form a proper class and therefore cannot be used to \Rcomment{fully represent the} ordinals in realizability models. Finally we present some tools for preserving cardinals in realizability models, including an analogue for realizability algebras of the forcing property known as the $\kappa$-chain condition. 
\end{abstract}

\section{Introduction}

\footnotetext{\textbf{Key words:} Logic, Classical Realizability, Set Theory}
\footnotetext{\Lcomment{\textbf{Acknowledgements}
This research was supported by RFSI (Réseau Francilien en Sciences Informatiques) and by ANR (Agence Nationale de la Recherche).  
}}

Realizability aims to extract the computational content of mathematical theories: a theory (or a logical system) is interpreted in a model of computation by establishing a correspondence between formulae and programs in a way that is compatible with the rules of deduction. For instance, a realizer of an implication $A\imp B$ is a program which, \Rcomment{when} applied to a realizer of $A$, returns a realizer of $B$. 
The origins of realizability \Rcomment{date} back to Kleene’s work in constructive mathematics \cite{Kleene1945}: Kleene’s realizability formalized the intuitionistic view that proofs are algorithms (computable functions) by interpreting proofs in Heyting arithmetic as recursive functions. 

In modern realizability, recursive functions are replaced by programs formalised in some variant of $\lambda$-calculus \cite{Barendregt}\Rcomment{. Modern} realizability \Lcomment{extends} the Curry-Howard isomorphism\Rcomment{, which is a method to establish} an isomorphism between proofs in intuitionistic logic and the terms of the simply typed $\lambda$-calculus. In the 90’s, the work of T. Griffin \cite{Griffin1989} led to pass the barrier of intuitionistic logic and to extend the Curry-Howard correspondence to classical logic by using the $\lambda_c$-calculus, an extension of $\lambda$-calculus that formalizes computation in the programming language Scheme. Building on the work of Griffin, the French mathematician J.-L. Krivine developed a method for realizing not only classical logic, but even Zermelo-Fraenkel set theory ZF (see for instance \cite{Krivine2001} and \cite{Krivine2012}). Following this technique it is possible to extract from any proof in ZF a program formalized in some variant of $\lambda$-calculus: for instance, the axiom of foundation is realized by the Turing fixed point combinator.
Krivine's technique generalizes the method of Forcing, but forcing models are uninteresting from a computational perspective since everything is realized by a single program. In this paper when we refer to realizability models, we will refer to realizability models not built from forcing unless explicitly specified.

One of the main obstacles when working with realizability models for set theory is the analysis of ordinals, in fact when one works in forcing (at least with posets), the resulting model is a proper extension of the ground model that shares \emph{the same ordinals}, this is not necessarily the case for realizability models. 

When working with forcing, \Rcomment{it is important to consider what happens to the cardinality of each ordinal in the forcing extension. For} instance, uncountable ordinals (hence uncountable cardinals) may collapse to countable ones\Rcomment{, a useful technique for producing many} different consistency results. Collapsing cardinals may happen as well in realizability models for set theory\Rcomment{, for example see Section 3 of \cite{Krivine2018}. In forcing, the simplest way to prevent these collapses is through \emph{chain conditions}, which we will generalise to realizability models. Finally,} ordinals code well-ordered sets, thus being able to talk about ordinals is especially important for realizing fragments of the axiom of choice. 

In this paper we present several results about the ordinals of realizability model\Rcomment{s} for set theory.

\subsection{Introduction to Realizability}

In this section we present \emph{Realizability Algebras}, which are the main building blocks for the construction of realizability models for set theory. We shall briefly explain the main intuition behind this construction. We start with a model of set theory, and we will use programs and stacks to evaluate the potential truth and falsity values of set theoretic statements. For computational reasons, we work with a non-extensional version of set theory, called $\ZFepsilon$, that involves two membership relations: the usual one and a \Rcomment{stricter} relation \Lcomment{that does not satisfies the axiom of extensionality}. We will use the terms of the $\lambda_c$-calculus, a variant of $\lambda$-calculus that \Rcomment{includes as a term the control operator \emph{call-with-current-continuation}; this additional term will be essential in showing that our resulting models satisfy a classical theory.} \Lcomment{We will use $\lambda_c$-terms} to evaluate the truth value\Rcomment{s} of formulas in the language of $\ZFepsilon$, \Lcomment{and we will use stacks, namely sequences of $\lambda_c$-terms, to evaluate the falsity values of such formulas.} \Rcomment{While related to each other, truth values and falsity values will be} \Lcomment{fundamentally different,} so that a $\lambda_c$-term is in the truth value of a formula (we say that it ``realizes the formula''), if it is somehow ``incompatible'' with every stack in the falsity value of that formula\Lcomment{, but in general the falsity value of a formula $\varphi$ will not necessarily coincide with the truth values of $\lnot \varphi.$} These definitions will respect certain logical constraints such as no stack can be in the falsity value of $\top,$ and every stack is in the falsity value of $\perp$. 
Then we choose some privileged collection of $\lambda_c$-terms that we call \emph{realizers}, and we will show that the set of formulas that are realized by some realizer forms a consistent theory which includes $\ZFepsilon$ and is closed under the rules of derivation of classical natural deduction. Finally, a realizability model will be a model of such a theory. Since $\ZFepsilon$ is a conservative extension of \tf{ZF}, such a model will induce a model of \tf{ZF}.\\ 

The main ingredients of realizability algebras are \emph{$\lambda_c$-terms}, \emph{stacks} and \emph{processes} which we define next. We will give the definition in full generality, in particular allowing for non-empty sets of \emph{special instructions}. These are customisable constants which can be added to our realizability algebras to ensure the models satisfy additional principles. For example, if the algebra is countable and contains the special instruction \emph{quote} then one can prove Dependent Choice holds in the model. \Rcomment{It is worthwhile to note that quote only makes sense if the algebra is countable. Therefore, the inclusion of additional special instructions can significantly change the theory we are able to realize.} However, \Rcomment{except for \Cref{Section:ChainConditions},} all the statements in this paper will be realized by terms of the $\lambda_c$-calculus without any special instruction\Rcomment{. From this, it will follow that they hold in every realizability model, even if the calculus includes special instructions.}

\begin{definition} Let \tf{V} be a model of \tf{ZF} and let $A, B$ be two sets in \tf{V} \Lcomment{with $B\neq \emptyset$:}

\begin{itemize}
\item We let $\Lambda^{\text{open}}_{A,B}$ and $\Pi_{A,B}$ denote the elements of \tf{V} defined by the following grammars, modulo $\alpha$-equivalence. Their elements are called respectively \emph{$\lambda_c$-terms} and \emph{stacks}:
 
$\begin{array}{rrlll}

\multicolumn{5}{l}{\Lambda^{\text{open}}_{A,B}~ (\textrm{\bf $\lambda_c$-terms}):}\\
& t, s\ ::=&  & x & \textrm{(variable; we choose a set of variables that is countable in \tf{V})} \\
&&\vert & ts & \textrm{(application)}\\
&&\vert &  \lambda u \lambdaapp t & \textrm{(abstraction; $u$ is a variable and $t$ is a $\lambda_c$-term)}\\
&&\vert & \cc & \textrm{(call-with-current-continuation)}\\
&&\vert & \saverlz{\pi} & \textrm{(continuation constants; $\pi$ is a stack)}\\
&&\vert  & \xi_\alpha & \textrm{(special instructions; $\alpha \in A$)}
\end{array}$

$\begin{array}{rrlll}
\multicolumn{5}{l}{\Pi_{A,B}~ (\textrm{\bf Stacks}):}\\
 & \pi\ ::=  & & \omega_\beta & \textrm{(stack bottoms; $\beta \in B$)}\\ 
&&\vert & t \stackapp \pi & \textrm{($t$ is a closed $\lambda_c$-term and $\pi$ is a stack)} 
\end{array}$

 \item $\Lambda_{A,B} \in \tf{V}$ denotes the set of all closed $\lambda_c$-terms,
 \item $\mathcal{R}_{A,B} \in \tf{V}$ denotes the set of all closed $\lambda_c$-terms that contain no occurrence of a continuation constant. Such terms are called \emph{realizers}.
 
\item $\Lambda_{A,B} \star \Pi_{A,B} \in \tf{V}$ denotes the Cartesian product $\Lambda_{A,B} \times \Pi_{A,B}$. Its elements are called \emph{processes}. We will write $t \star \pi$ for $(t, \pi) \in \Lambda_{A, B} \star \Pi_{A, B}$.

\end{itemize}
\end{definition}

Application on $\lambda_c$-terms is left associative\Rcomment{, so} $(ts_1s_2\cdots s_n)$ means $\Rcomment{(}(\cdots ((ts_1)s_2)\cdots )s_n)$ and has higher priority than abstraction\Rcomment{.} ($\lamAbstOne ts)$ means $\Rcomment{(}\lamAbstOne (ts))$. When there is no ambiguity, we will drop the indexes A, B and simply denote \Rcomment{$\Lambda_{A,B}$} by $\Lambda,$ \Lcomment{$\Lambda_{A,B}^{open}$ by $\Lambda^{open}$,} $\Pi_{A,B}$ by $\Pi$ \Lcomment{and $\mathcal{R}_{A,B}$ by $\mathcal{R}$}.\\

The terms \Lcomment{$\cc$ and $\saverlz{\pi}$ refer to the} \emph{call-with-current-continuation} \Lcomment{procedure from the Scheme programming language, that allows access to the current context of a process (i.e. the current stack). Roughly speaking, $\cc$ saves the current context as a continuation term $\saverlz{\pi}$ that can be later restored; when a continuation term $\saverlz{\pi}$ is applied to a term, the existing context is eliminated and the applied continuation context is restored in its place, so that the process flow will resume at the point at which the context was captured. This behaviour is formalized in the definition of the }\emph{evaluation}\Lcomment{ hereafter, that sets the }rules of reduction on the set of processes. 

\begin{definition}
Let \tf{V} be a model of \tf{ZF} and let $A, B$ be two sets in \tf{V}. 
\begin{itemize}
\item $\prec_{A,B} \; \in \tf{V}$ is called the \emph{evaluation preorder} and denotes the smallest preorder on \hbox{$\Lambda_{A,B} \star \Pi_{A,B}$} such that:
\[
{ \arraycolsep=1.5pt \begin{array}{rclcrcll}
ts & \star & \pi & \quad \succ_{A,B} \quad & t & \star & s \point \pi \quad &\textrm{ (push)}\\
\lambda u \lambdaapp t & \star & s \point \pi & \quad \succ_{A,B} \quad & \quad  t[u \coloneqq s] & \star & \pi &\textrm{ (grab)}\\
\cc & \star & t\point \pi & \quad \succ_{A,B} \quad & t & \star & \saverlz{\pi} \point \pi &\textrm{ (save)}\\
\saverlz{\sigma} & \star & t\point \pi & \quad \succ_{A,B} \quad & t & \star & {\sigma} &\textrm{ (restore).}
\end{array} }
\]
\end{itemize}
\end{definition}

Note that there is no evaluation rule for the special instructions, depending on the context we may define other evaluation relations with specific evaluation rules for the special instructions. 
When there is no ambiguity, we will drop the indexes $A, B$ and simply write \Rcomment{$\prec$} for \Rcomment{$\prec_{A, B}$}.

\begin{definition}
Let \tf{V} be a model of \tf{ZF}. A \emph{realizability algebra} in \tf{V} is a tuple $\mathcal{A}=(A, B, \Rcomment{\prec_{A,B}, \;} \Perp)$ such that:  
 \begin{itemize} 
\item $\mathcal{A} \in \tf{V}$ (i.e.\ $A \in \tf{V}$, $B \in \tf{V}$, $\Rcomment{\prec_{A,B} \, \in \tf{V}}$ and $\Perp \in \tf{V}$);
\Rcomment{\item $\Rcomment{\prec_{A,B}}$ is an (possibly trivial) extension of the evaluation preorder on $\Lambda_{A,B} \star \Pi_{A,B}$.}
\item $\Perp$ is a subset of $\Lambda_{A,B} \star \Pi_{A,B}$ that is \Lcomment{upwards closed} for $\succ_{A,B}$, i.e. if $t\star \pi\succ_{A,B} s \star \sigma$ and $s \star \sigma \in \Perp$, then $t\star \pi\in \Perp.$ It is called the \emph{pole} of the realizability algebra. 
 \end{itemize}
\end{definition}

\Lcomment{The pole is not in general downward closed. A downward closed pole induces the following property which is true in Kleene's realizability but not necessarily in classical realizability: if $tu$ realizes $B$ for every $u$ that realizes $A,$ then $t$ realizes $A\imp B.$}



\medskip

In order to define a realizability model for classical set theory, we consider a non-extensional conservative extension of the usual set theory. This theory was originally formulated by Friedman in \cite{Friedman1973} in his proof that \tf{ZF} is equiconsistent with \tf{IZF} \Rcomment{using a double negation translation. However, 
this translation does not work well with the Axiom of Extensionality. For this reason, the proof first goes via a non-extensional set theory which} contains two distinct membership relations: $\in$ which behaves like the standard membership relation, and $\rlzin$ which is a form of ``\emph{strong membership}''.

Throughout this paper, we will work in first-order logic \emph{without equality}: \Rcomment{individual} language may contain a symbol that happens to be written ``\(=\)'', but 
models are not required to interpret it by ``meta'' equality. In addition, we will assume that the only primitive logical constructions are $\rightarrow$, $\top$, $\bot$, and $\forall$; for $\vee$, $\wedge$, and $\exists$, we will use De Morgan's encoding. \Rcomment{We remark here that logical connectives will be evaluated in their usual, right associate, manner.} Thus:
\begin{itemize}
\item $\varphi \wedge \psi$ means $(\varphi \rightarrow \Rcomment{(} \psi \rightarrow \perp \Rcomment{)}) \rightarrow \perp$,
\item $\varphi \vee \psi$ means $(\varphi \rightarrow \perp) \rightarrow \Rcomment{(}(\psi \rightarrow \perp) \rightarrow \perp\Rcomment{)}$,
\item $\exists x~ \varphi(x)$ means $(\forall x~ (\varphi(x) \rightarrow \perp)) \rightarrow \perp$.
\end{itemize}

We will denote by $\mathcal{L}_{\in}$ the first-order language over the signature \(\{\in, \simeq\}\) where \(\in\) and \(\simeq\) are binary relation symbols. The language of $\ZFepsilon$ requires two distinct symbols for the membership relation, $\in$ and $\rlzin$ (the former will refer to the usual extensional membership relation, the latter will correspond to a \Rcomment{stricter} non-extensional membership relation). However, for computational reasons it is better to take as primitives the negative versions of those symbols. \Lcomment{In fact we are going to define truth and falsity values of formulas, intuitively the falsity values can be regarded as computations that challenge the validity of the formula, but to compute the falsity of an expression of the form $a\notin b$ is easier than to compute the falsity of $a\in b,$ since for the former it is enough to show that $a$ is in fact in $b,$ while for the latter one would need to go over all the elements of $b$ and check that $a$ is none of them, this computation may not terminate.} 
Thus the language of $\ZFepsilon$, which is denoted $\mathcal{L}_{\rlzin},$ is the first-order language over the signature \(\{\notin, \subseteq, \notrlzin, \neq\}\), where all \(4\) symbols are binary relation symbols. It can be proven that $\neq$ is definable from $\notrlzin$ via the Leibniz equality and is therefore not necessary in the signature, however we include it here for practical purposes. $\mathrm{Fml}_\in$ and $\mathrm{Fml}_{\rlzin}$ denote the collection of all formulas in $\mathcal{L}_\in$ and $\mathcal{L}_{\rlzin}$ respectively. In the language $\mathcal{L}_{\rlzin}$, we will use the following abbreviations:

\begin{center}
\begin{tabular}{| r | l || r | l |}
\hline
Abbreviation & Meaning & Abbreviation & Meaning\\  
\hline
$a \rlzin b$ & $a \notrlzin b \rightarrow \perp$  & $a \simeq b$ & $(a \subseteq b) \wedge (b \subseteq a)$\\
$a \in b$ & $a \not\in b \rightarrow \perp$ & $\forall x \rlzin a~ \varphi(x)$ & $\forall x~ (x \rlzin a \rightarrow \varphi(x))$\\
$a = b$ & $a \neq b \rightarrow \perp$ & $\exists x \rlzin a~ \varphi(x)$ & $(\forall x~ (\varphi(x) \rightarrow x \notrlzin a)) \rightarrow \bot$\\
\hline
\end{tabular}
\end{center}

In particular, by a slight abuse of notation, we will consider $\mathrm{Fml}_{\in}$ to be a subset of $\mathrm{Fml}_{\rlzin}$. \(\tf{ZF}\) denotes the usual set theory, written in the language $\mathcal{L}_\in$, i.e.\ \(\tf{ZF}\) is a subset of $\mathrm{Fml}_\in,$ while $\ZFepsilon$ denotes non-extensional set theory, as defined by Krivine, written in the language $\mathcal{L}_{\rlzin}$ (i.e.\ $\ZFepsilon$ is a subset of $\mathrm{Fml}_{\rlzin}$). In a nutshell, the axioms of $\ZFepsilon$ state that:
\begin{itemize}
\item An equivalent presentation of the axioms of \tf{ZF} \emph{minus the Axiom of Extensionality} (essentially the double negation) are satisfied \emph{over the signature} $\{ \notrlzin, \neq \}$ (rather than $\{\in, \simeq \}$).
\item \(\in\) is the extensional collapse of \(\rlzin\): $x\in y$  \textit{iff} there is $x'\rlzin y$ such that $x\simeq x'$; 
\item $\subseteq$ is the extensional inclusion: $x\subseteq y$  \textit{iff} for every $z\rlzin x,$ we have $z\in y$;
\item \(\simeq\) is extensional \Rcomment{equality}: two sets are $\simeq$-equal \textit{iff} they have the same $\in$-elements.
\end{itemize} 
For full details, including the list of the axioms of $\mathcal{L}_{\rlzin}$, we refer the reader to \cite{Krivine2012}; see also Friedman's earlier account in \cite{Friedman1973}.
As proven in \cite{Friedman1973}, $\ZFepsilon$ is a conservative extension of $\tf{ZF}$:

\begin{theorem}[Friedman / Krivine] \label{theorem:ZFepsilonToZF}
Let $\varphi$ be a closed formula in $\mathcal{L}_\in$: $\varphi$ is a consequence of \(\tf{ZF}\) if and only if it is a consequence of $\ZFepsilon$.
\end{theorem}

A proof of this fact can be found in \cite{Krivine2001}. For further details we refer to \cite{Matthews2023}. 

Whenever $\mathcal{L}$ is a first-order language that contains $\mathcal{L}_{\rlzin}$, we will denote by $\ZFepsilonL{\mathcal{L}}$ the theory obtained by taking $\ZFepsilon$ and enriching all the axiom schemas to include the formulas of $\mathcal{L}$.

Our construction of realizability models follows the presentation in \cite{Matthews2023}. Let \tf{V} be a model of Zermelo-Fr\ae{}nkel set theory, \tf{ZF}, and let $\mathcal{A} = \Rcomment{(A, B, \prec, \Perp)}$ be a realizability algebra in $\tf{V}.$ 

We define $\rlzstr^{\mathcal{A}, \tf{V}} \subseteq \tf{V}$ as follows: 
for any ordinal $\alpha \in \tf{V}$, let $\rlzstr_\alpha^{\mathcal{A}, \tf{V}} \coloneqq \bigcup_{\beta < \alpha} \mathcal{P}(\rlzstr_\beta^{\mathcal{A}, \tf{V}} \times \Pi),$  then let $\rlzstr^{\mathcal{A}, \tf{V}} \coloneqq \bigcup_{\alpha\in \text{Ord}} \rlzstr_\alpha^{\mathcal{A}, \tf{V}},$ where $\text{Ord}$ denotes the class of ordinals in \tf{V}. The elements of $\rlzstr^{\mathcal{A}, \tf{V}}$ are called $(\mathcal A, \tf{V})$-\emph{names}. Note that for all $\alpha$, $\rlzstr_\alpha^{\mathcal{A}, \tf{V}} \in \tf{V}$, but $\rlzstr^{\mathcal{A}, \tf{V}} \notin \tf{V}$ ($\rlzstr^{\mathcal{A}, \tf{V}}$ is a \emph{proper class} in \tf{V}). We will generally drop the exponents and simply write $\rlzstr_\alpha$ and $\rlzstr$. Given an element $a \in \rlzstr$ we let $\ff{dom}(a) \coloneqq \{ b \divline \exists \pi \in \Pi ~ (b, \pi) \in a \} \in \tf{V}$.

Certain names will play a special role in the realizability model: the so-called \emph{gimel names} denoted \Rcomment{$\cjgimel(a)$} and \emph{reish names} denoted $\fullname a$ \Rcomment{are} defined as follows.

\begin{definition}
For $x \subseteq \rlzstr$, we define $\cjgimel(x):= x \times \Pi.$
\end{definition}

\begin{definition}
    We define $\fullname{x}$ for $x \in \tf{V}$ recursively as:
    \[
    \fullname{x} \coloneqq \{ (\fullname{y}, \pi) \divline y \in x\Rcomment{, \pi \in \Pi} \}. 
    \]
\end{definition}

A function $f : \rlzstr^n \to \rlzstr$ is said to be \emph{$\mathcal A$-definable} if there is a formula $\varphi(x_1,\ldots,x_n, y) \in \mathrm{Fml}_{\rlzin}$ 
such that, for any $a_1, \ldots, a_n \in \rlzstr$ and $b \in \tf{V}$, $\tf{V} \models \varphi(a_1,\ldots,a_n,b)$ if and only if $b = f(a_1, \ldots, a_n)$. Let $\mathcal{L}_{\rlzin}^{\mathcal A}$ be the language obtained from $\mathcal{L}_{\rlzin}$ by adding for each ${\mathcal A}$-definable function $f : \rlzstr^n \to \rlzstr$ an $n$-ary function symbol ``$f$''.

The \emph{realizability interpretation} of $\mathcal{L}_{\rlzin}^{\mathcal A}$ in $\mathcal{A}$ consists of the following. 

\begin{definition}\label{def:truth-falsity-values}

To each closed formula $\varphi$ in $\mathcal{L}_{\rlzin}^{\mathcal A}$ with parameters in $\rlzstr$, we associate a \emph{truth value} $\verity{\varphi} \subseteq \Lambda$ and a \emph{falsity value} $\falsity{\varphi} \subseteq \Pi,$ they are defined jointly by induction on the complexity of $\varphi$:


    \begin{itemize}
    \item $\verity{\varphi} \coloneqq \{ t \in \Lambda \divline \forall \pi \in \falsity{\varphi}, t \star \pi \in \Perp \}$;
    \item $\falsity{\top} \coloneqq \emptyset$ and $\falsity{\perp} \coloneqq \Pi$;
    \item $\falsity{a \notrlzin b} \coloneqq \{ \pi \in \Pi \divline (a, \pi) \in b \};$ 
    \item $\falsity{a \neq b} := \falsity{\top}$ if $a\neq b$, $\falsity{\bot}$ otherwise;
    \item $\falsity{a \not\in b} \coloneqq \displaystyle{\bigcup}_{c \in \ff{dom}(b)} \{ t \point t' \point \pi \divline (c, \pi) \in b, t \in \verity{a \subseteq c}, t' \in \verity{c \subseteq a} \}$;
    \item $\falsity{a \subseteq b} \coloneqq \displaystyle{\bigcup}_{c \in \ff{dom}(a)} \{ t \point \pi \divline (c, \pi) \in a, t \in \verity{c \not\in b} \}$;
    \item $\falsity{\psi \rightarrow \theta} \coloneqq \{ t \point \pi \divline t \in \verity{\psi}, \pi \in \falsity{\theta} \}$;
    \item $\falsity{\forall x~ \varphi(x)} \coloneqq {\displaystyle \bigcup}_{a \in \rlzstr} \falsity{\varphi[a / x]}$.
\end{itemize}

Formally, $\falsity{a \not\in b}$ and $\falsity{a \subseteq b}$ are defined by induction on the pair $(\max(\ff{rk}_\rlzstr(a), \ff{rk}_\rlzstr(b)),\allowbreak \min(\ff{rk}_\rlzstr(a), \ff{rk}_\rlzstr(b)))$ under the product order, where $\ff{rk}_\rlzstr(c) := \min\,\{\alpha \divline c \in \rlzstr_{\alpha+1}\}$. 

We say that a closed $\lambda_c$-term $t$ \emph{realizes} a closed formula $\varphi$ with parameters in $\rlzstr$
and write $t \Vdash \varphi,$ whenever $t \in \verity{\varphi}$.
\end{definition}

\Rcomment{As a first example of using realizers, we shall show how to realize \emph{Peirce's Law}. As this is logically equivalent to the law of Excluded Middle, from this it will follow that we will obtain a classical theory.

\begin{proposition} \label{theorem:KpiAndNegation}
    Suppose that $\pi \in \falsity{\varphi}$. Then for any formula $\psi$ in $\mathcal{L}_{\rlzin}^{\mathcal A}$, $\saverlz{\pi} \Vdash \varphi \rightarrow \psi$. In particular, $\saverlz{\pi} \Vdash \neg \varphi$.
\end{proposition}

\begin{proof}
    Suppose that $t \Vdash \varphi$ and $\pi \in \falsity{\psi}$. Then, $\saverlz{\pi} \star t \stackapp \pi \succ t \star \pi$. Thus, since $t \Vdash \varphi$ implies that $t \star \pi \in \Perp$, $\saverlz{\pi} \star t \stackapp \pi \in \Perp$. 
\end{proof}

\begin{proposition}[Peirce's Law]
    For any formulas $\varphi$ and $\psi$ in $\mathcal{L}_{\rlzin}^{\mathcal A}$, $\cc \Vdash ((\varphi \rightarrow \psi) \rightarrow \varphi) \rightarrow \varphi.$
\end{proposition}

\begin{proof}
    Suppose that $t \Vdash (\varphi \rightarrow \psi) \rightarrow \varphi$ and $\pi \in \falsity{\varphi}.$ By \Cref{theorem:KpiAndNegation}, $\saverlz{\pi} \Vdash \varphi \rightarrow \psi$, from which it follows that 
    \[
    \cc \star t \stackapp \pi \succ t \star \saverlz{\pi} \stackapp \pi \in \Perp.
    \]
\end{proof}}

Now, we would like to associate to $\mathcal A$ a ``realizability theory'' consisting of all closed formulas which are realized. However, \Rcomment{by the above argument}, for all $t \star \pi \in \Perp$ the $\lambda_c$-term $\saverlz{\pi}\,t$ realizes the formula $\bot$. Therefore, in order to obtain a realizability theory that is not automatically inconsistent, we will need to exclude terms of this shape; this is where the set $\mathcal{R}$ of \emph{realizers} comes into play (i.e. the closed $\lambda_c$ terms containing no continuation constant): 

\begin{definition}
The \emph{realizability theory of $(\mathcal{A}, \tf{V})$}, denoted by $T_{\mathcal{A},\tf{V}}$, is the set of all closed formulas, $\varphi$, of $\mathcal{L}_{\rlzin}^{\mathcal{A}}$ with parameters in $\rlzstr$ \Rcomment{for which} there exists $t \in \mathcal{R}$ such that $t$ realizes $\varphi$.
\end{definition}

The following facts are standard (see e.g.\ \cite{Krivine2012}):\begin{itemize}
\item the realizability theory of $(\mathcal{A}, \tf{V})$ is closed under classical deduction, (i.e. if $\varphi \in T_{\mathcal{A},\tf{V}}$ and $\varphi$ entails $\psi$ in classical logic, then $\psi \in T_{\mathcal{A},\tf{V}}$);
\item this theory is consistent if and only if for every $t\in\mathcal{R}$ there is a stack $\pi$ such that $t \star \pi \not\in \Perp$;
\item this theory is generally not complete.
\end{itemize}



\begin{theorem}[Krivine]
    Let \tf{V} be a model of \tf{ZF} and $\mathcal A$ a realizability algebra in \tf{V}. The realizability theory of $(\mathcal{A}, \tf{V})$ contains $\ZFepsilonL{\mathcal{L}_{\rlzin}^{\mathcal{A}}}$. In particular, it contains $\ZFepsilon$, and therefore \tf{ZF}.
\end{theorem}

We refer to \cite{Krivine2012} for a proof of this, or \cite{Matthews2023} for an alternative proof using the setup given. This justifies the following definition:
\begin{definition}
A \emph{realizability model of $\ZFepsilon$} is a pair $\rlzmodel = (\tf{V}, \mathcal A)$, with \tf{V} a model of \tf{ZF} and $\mathcal A$ a realizability algebra in \tf{V}. We write $\rlzmodel \Vdash \varphi$ for $``T_{\mathcal{A},\tf{V}}$ contains $\varphi$''.
\end{definition}

Sometimes, we will argue within models of the realizability theory $T_{\mathcal{A},\tf{V}}$ and by abuse of language we will call \emph{the} realizability model any model of $T_{\mathcal{A},\tf{V}}.$ \Rcomment{As previously mentioned, the realizability models derived from forcing posets are computationally uninteresting since everything is realized by a single program. From this it can be proven that the set of truth values, $\fullname{2}$, consists of precisely two $\rlzin$-elements. It can also be proven that if $\mathcal{A}$ is a countable realizability algebra and $\fullname{2}$ consists of precisely two $\rlzin$-elements then the resulting model is equivalent to a forcing model. We refer to \cite{Matthews2023} for full details. We shall therefore say that either a realizability model or $\fullname{2}$ is \emph{trivial} whenever $\fullname{2}$ consists of precisely two $\rlzin$-elements.}


We end this introduction with two brief pieces of notation that will be frequently used throughout this work\Rcomment{, the first of which is a variant of subsets for just the $\rlzin$-structure.} 

\begin{definition}
    Given $a,b\in N,$ we let $a \subseteq_{\rlzin} b \equiv \forall x (x \rlzin a \rightarrow x \rlzin b)$.
\end{definition}

\Rcomment{The second piece of notation will give a significantly easier way to view bounded universal quantifiers. We will introduce a new formula $\forall x^{\cjgimel{A}} \varphi(x)$ which will be used to realize formulas of the form $\forall x \rlzin \cjgimel{A} \varphi(x)$ whenever $A$ is a collection of names. In particular, this will work for $\fullname{a}$ for every $a \in \rlzstr$.}

\begin{definition} \label{definition:BoundedUniversals}
    Given $\Rcomment{A} \subseteq N$ we define $\falsity{\forall x^{\Rcomment{\cjgimel(A)}} \varphi(x)} =  \bigcup_{b \in \Rcomment{A}} \falsity{\varphi(b)}$.
\end{definition}

\Rcomment{The benefit of this presentation for bounded universal quantifiers is that it behaves in an analogous way to its unbounded counterpart. Namely, it is formed as the union of the falsifier of the formula for each element in $A$. Moreover, using double negations, the formula $\forall x \rlzin \cjgimel(A) \varphi(x)$ is in fact $\forall x (\neg \varphi(x) \rightarrow x \notrlzin \cjgimel(A))$ which can lead to much less readable sentences. We remark here that we required our names to be of the form $\cjgimel(A)$ rather than arbitrary names in order for \cref{theorem:BoundedUniversals:StateToAbbr} of \Cref{theorem:BoundedUniversals} to go through.}

\begin{proposition}  \label{theorem:BoundedUniversals} \Rcomment{Fix $A \subseteq \rlzstr$.}
    \begin{enumerate}
        \item $\lamAbstOne \lamAbstTwo \app{\lamtermTwo}{\lamtermOne} \Vdash \forall x^{\Rcomment{\cjgimel(A)}} \varphi(x) \rightarrow \forall x (\neg \varphi(x) \rightarrow x \notrlzin \Rcomment{\cjgimel(A)})$.
        \item \label{theorem:BoundedUniversals:StateToAbbr}
        $\lamAbstOne \rapp{\cc}{\lamAbst{k} \app{u}{k}} \Vdash \forall x (\neg \varphi(x) \rightarrow x \notrlzin \Rcomment{\cjgimel(A)}) \rightarrow \forall x^{\Rcomment{\cjgimel(A)}} \varphi(x)$.
        \item \Rcomment{Hence, for every set $c \in \tf{V}$, $\rlzmodel \Vdash \forall x^{\fullname{c}} \varphi(x) \leftrightarrow \forall x \rlzin \fullname{c} \varphi(x)$.}
    \end{enumerate}
\end{proposition}

\begin{proof} The first \Rcomment{and third claims are straight forward}.

For the second claim, let $t\Vdash \forall x (\neg \varphi(x) \rightarrow x \notrlzin \Rcomment{\cjgimel{A}})$ and $\pi\in \falsity{\varphi(b)}$ for some $b\in \Rcomment{A}$. \Rcomment{Since $(b, \sigma) \in \cjgimel(A)$ for every $\sigma \in \Pi$,} $\pi\in \falsity{b\notrlzin a}.$ It follows that $k_\pi\Vdash \lnot \varphi(b)$ and $t\ast k_\pi\stackapp \pi\in 
\Perp.$ \Rcomment{Hence we} have $\lamAbstOne \rapp{\cc}{\lamAbst{k} \app{u}{k}}\ast t\stackapp \pi\succ t\ast k_\pi\stackapp \pi\in \Perp$ from which we can conclude that $\lamAbstOne \rapp{\cc}{\lamAbst{k} \app{u}{k}} \Vdash \forall x (\neg \varphi(x) \rightarrow x \notrlzin {\Rcomment{\cjgimel(A)}}) \rightarrow \forall x^{\Rcomment{\cjgimel(A)}} \varphi(x).$  
\end{proof}

\subsection{Useful Realizers}

Here we briefly list some of the realizers that we will make frequent use of throughout this paper. For some of the simpler statements we will also explicitly state $\lambda$-terms which satisfy the required properties and leave it to the interested reader to verify this. However, we remark here that all of these results follow from the fact that the realizability model satisfies $\ZFepsilon$.

\begin{itemize}
    \item $\identity = \lamAbstOne \lamtermOne$.
    \item Let $\rlzXSubsetX = \app{\theta}{\theta}$, where $\theta = \lamAbstOne \lamAbstTwo \twoapp{\rapp{\lamtermTwo}{\app{\lamtermOne}{\lamtermOne}}}{\app{\lamtermOne}{\lamtermOne}}$.\footnote{This can be seen as a type of ``Turing fixed point combinator''.} Then $\rlzXSubsetX \Vdash \forall x (x \subseteq x)$.
    \item Let $\rlzXSimeqX = \lamAbstOne \lapp{\app{\lamtermOne}{\rlzXSubsetX}}{\rlzXSubsetX}$. Then $\rlzXSimeqX \Vdash \forall x (x \simeq x)$, where $x \simeq y$ is an abbreviation for \linebreak[4] $(x \subseteq y \rightarrow \Rcomment{(} y \subseteq x \rightarrow \perp \Rcomment{)}) \rightarrow \perp$.
    \item Let $\rlzXNotInX = \app{\theta'}{\theta'}$, where $\theta' = \lamAbstOne \lamAbstTwo \lamAbstThree \rapp{\lamtermTwo}{\app{\lamtermOne}{\lamtermOne}}$. Then    $\rlzXNotInX \Vdash \forall x (x \not\in x)$.
    \item $ \rlzSubsetTrans \Vdash \forall x \forall y \forall z (x \subseteq y \rightarrow (y \subseteq z \rightarrow x \subseteq z))$.
    \item $\rlzInIsNotSimeq \Vdash \forall x \forall y (x \rlzin y \rightarrow x \not\simeq y)$.
    \item Let $\rlzNotNot = \lamAbstOne \lamAbstTwo \app{\lamtermTwo}{\lamtermOne}$. Then $\rlzNotNot \Vdash \varphi \rightarrow ((\varphi \rightarrow \perp) \rightarrow \perp)$.
\end{itemize}

\Rcomment{We also have a standard way to encode the natural numbers in $\lambda$-calculus, known as the \emph{Church numerals}}.  

\begin{definition}
\Lcomment{(Church numerals) Let $\underline{0}$ denote the $\lambda$-term $\lamAbstOne \lamAbstTwo \lamtermTwo.$ For $n \in \omega$, $\underline{n+1}$ denotes the $\lambda$-term $\lamAbstOne \lamAbstTwo \twoapp{\app{\underline{n}{\lamtermOne}}}{\app{\lamtermOne}{\lamtermTwo}}$. In particular, 
    $\underline{1}$ denotes the $\lambda$-term $\lamAbstOne \lamAbstTwo \app{\lamtermOne}{\lamtermTwo}$.}
\end{definition}
\Rcomment{
\begin{remark}
    Let $s = \lamAbst{n} \lamAbstOne \lamAbstTwo \twoapp{\app{n}{\lamtermOne}}{\app{\lamtermOne}{\lamtermTwo}}$. Then, for every $n \in \omega$ we have $s\underline{n} \rightarrow_\beta \underline{n+1}$.
\end{remark}}

\section{Functions} \label{Section:Functions}

In the latter half of this paper we will discuss cardinality of objects in the realizability model. In order to do this, we need to have a clear concept of what a function is. This has different consequences depending on where extensional and non-extensional terms are used. Thus we must be very careful to differentiate between each version, where $\in$- will denote the standard definition over \tf{ZF}. \Rcomment{Each of these definitions can also be generalised in the obvious way to $n$-ary functions.}

\begin{definition}
    ($\ZFepsilon$) $f$ is said to be an $\in$-function with domain $a$ if it is a binary relation satisfying the following two conditions:
    \vspace{-5pt}
    \begin{itemize} \setlength \itemsep{2pt}
        \item ($\in$-totality) $\forall x \in a \, \exists y \, f(x,y);$
        \item ($\in\text{-}\ff{Fun}(f, a)$) $\forall x \in a \, \forall y, y' \big( f(x, y) \land f(\Rcomment{x}, y') \rightarrow y \simeq y' \big).$
    \end{itemize}
\end{definition}

\begin{definition}
    ($\ZFepsilon$) $f$ is said to be an $\rlzin$-function with domain $a$ if it is a binary relation satisfying the following two conditions:
    \vspace{-5pt}
    \begin{itemize} \setlength \itemsep{2pt}
        \item ($\rlzin$-totality) $\forall x \rlzin a \, \exists y \, f(x,y);$
        \item ($\rlzin\text{-}\ff{Fun}(f, a)$) $\forall x \rlzin a \, \forall y, y' \big( f(x, y) \land f(\Rcomment{x}, y') \rightarrow y = y' \big).$
    \end{itemize}
\end{definition}

\begin{definition}
($\ZFepsilon$) $f$ is said to be an \emph{extensional function with domain} $a$ if it is a binary relation satisfying the following two conditions:
\vspace{-5pt}
\begin{itemize} \setlength \itemsep{2pt}
    \item ($\rlzin$-totality) $\forall x \rlzin a \, \exists y \, f(x,y);$
    \item ($\ff{ExtFun}(f, a)$) $\forall x, x' \rlzin a \, \forall y, y' \big( f(x, y) \land f(x', y') \land x \simeq x' \rightarrow y \simeq y' \big).$
\end{itemize}
\end{definition}

\begin{definition}
\Rcomment{Let $f$ be a binary relation satisfying any of the above versions of being a function.} $f$ is said to be an \emph{extensional injection} if 
\[
\forall x, x' \rlzin a \, \forall y \, \big( f(x, y) \land f(x', y) \rightarrow x \simeq x' \big).
\]
$f$ is said to be an \emph{$\rlzin$-surjection} onto $b$ if
\[
\forall y \rlzin b \, \exists x \rlzin a \, f(x, y).
\]
\end{definition}

\begin{remark}
    Because \tf{ZF} is being interpreted in $\ZFepsilon$ with equality given by $\simeq$, when we have an extensional function $f(x, y)$ we will also read this as $f(x) \simeq y$.
\end{remark}

An important concept is how to build functions in a realizability model. In general, this is difficult to do for extensional functions. However, given a function $f$ in the ground model it is possible to ``\emph{lift}'' it in some precise sense to an $\rlzin$-function in the realizability model. We outline this method in this section, full details are then given in Section 16 of \cite{Matthews2023}.


We define special functions $\sing, \up, \op$ that in the realizability model will be interpreted as the singleton, the unordered pair, and the ordered pair respectively. 

\begin{definition} \label{definition:SngUpOp} \
\vspace{-5pt}
\begin{itemize} \setlength \itemsep{2pt}
        \item $\sing \colon \rlzstr \rightarrow \rlzstr$, $a \mapsto \{a\} \times \Rcomment{\Pi}$.
        \item $\up \colon \rlzstr \times \rlzstr \rightarrow \rlzstr$, $(a, b) \mapsto \{ (a, \underline{0} \stackapp \pi) \divline \pi \in \Pi \} \cup \{ (b, \underline{1} \stackapp \pi) \divline \pi \in \Pi \}$,
        \item $\op \colon \rlzstr \times \rlzstr \rightarrow \rlzstr$, $(a, b) \mapsto \up(\up(\sing(a), \fullname{0}), \sing(\sing(b)))$.
    \end{itemize}

\end{definition}

\begin{definition} \label{definition:FunctionLift}
    Suppose that $A \subseteq \rlzstr$ and $f \colon A \rightarrow \rlzstr$ is a \Rcomment{definable} function in \tf{V}. We define the \emph{lift} of $f$, $\lift{f}$ as
    \[
    \lift{f} \coloneqq \{ (\op(c, f(c)), \pi) \divline c \in A, \pi \in \Pi \}.
    \]
\end{definition}

\begin{proposition}[\cite{Matthews2023}]
    Let $A \subseteq \rlzstr$ and let $f \colon A \rightarrow \rlzstr$ be a \Rcomment{definable} function in \tf{V}. Then $\rlzmodel \Vdash \rlzin\text{-}\ff{Fun}(\lift{f}, \cjgimel(A))$ and for all $c \in A$,
    \[
    \Rcomment{\rlzmodel \Vdash \forall z (\op(c, z) \rlzin \lift{f} \longleftrightarrow z = f(c)).} 
    \]
\end{proposition}

\Rcomment{One particular benefit of this result is that it gives us the ability to give a canonical name for the basic internal operations in $\rlzmodel$. For example, if $A = \rlzstr$ and $f$ is a definable class function, such as the function $\sing$, then $\lift{f}$ will be a definable class $\rlzin$-function which outputs the desired operation. Using this, we can treat such functions as symbols that have been added to the language and freely use them in $\ZFepsilon$. The construction of $\lift{f}$ can also be generalised in a straightforward manner to functions with domain $A \times B$, where $A, B \subseteq \rlzstr$, from which we can lift the other two functions in \Cref{definition:SngUpOp}. We refer to Section 9 of \cite{Matthews2023} for full details.}

\begin{theorem}[\cite{Matthews2023}]
    The functions $\sing, \up, \op$ globally lift to the realizability model.
    Moreover, 
    \[
    \rlzmodel \Vdash \forall a, b (\sing(a) \simeq \{ a\} \land \up(a, b) \simeq \{a, b\} \land \op(a, b) \simeq \langle a, b \rangle),
    \]
    where $\langle a, b \rangle $ is the Wiener definition of the ordered pair, $\{ \{\{a\}, 0\}, \{\{b\}\} \}$.
\end{theorem}

Instead of having to work with ``\emph{the unique} $y$ \emph{such that} $\op(x, y) \rlzin \lift{f}$'' we \Rcomment{want to instead treat $\lift{f}$ as a function symbol. For this we introduce a new falsity value to work directly with $\lift{f}(x)$.}

\begin{definition}
    Suppose that $A \subseteq \rlzstr$ and $f \colon A \rightarrow \rlzstr$ is a function in \tf{V}. Let $\lift{f}$ be the lift of $f$. Then, for any formula $\varphi(u, v, \overrightarrow{w})$ and $\overrightarrow{w} \subset \rlzstr$,
    \[
    \falsity{\varphi(x, \lift{f}(x), \overrightarrow{w})} = \begin{cases}
        \falsity{\varphi(x, f(x), \overrightarrow{w})} & \text{ if } x \in A, \\
        \Rcomment{\Pi} & \text{ otherwise.}
    \end{cases}
    \]
\end{definition}

\begin{proposition}[\cite{Matthews2023}]
    For every formula $\varphi$,
    \[
    \rlzmodel \Vdash \forall x^{\cjgimel(A)} \big( \varphi(x, \lift{f}(x), \overrightarrow{w}) \longleftrightarrow \forall y ( \op(x, y) \rlzin \lift{f} \rightarrow \varphi(x, y, \overrightarrow{w})) \big).
    \]
\end{proposition}

The Axiom of Choice is hard to realize, instead a weaker version of the Axiom of Choice 
can be seen to hold in many realizability models, this is the so-called \emph{Non-extensional Axiom of Choice}, \tf{NEAC}


\begin{definition}
    \tf{NEAC} is the assertion that every $\rlzin$-relation can be reduced to an $\rlzin$-function. Namely: For every $r$, there exists some $f$ such that
    \begin{itemize}
        \item $\forall x \forall y, y' ( \op(x, y) \rlzin f \land \op(x, y') \rlzin f \rightarrow y = y')$,
        \item $f \subseteq_{\rlzin} r$,
        \item $\forall x \forall y \exists y' (\op(x, y) \rlzin r \rightarrow \op(x, y') \rlzin f)$.
    \end{itemize}
\end{definition}

As proven by Krivine (see for instance \cite{Krivine2012}), if the realizability algebra contains a special instruction $\rlzfont{q}$ implementing the instruction \emph{quote} \Rcomment{of} the LISP programming language, then the realizability model realizes \tf{NEAC}. 

\begin{lemma}
    Suppose that $\mathcal{A}$ is a countable realizability algebra and $\Lambda$ contains a special instruction $\rlzfont{q}$ satisfying the instruction \emph{quote}. Then $\rlzmodel \Vdash \tf{NEAC}$.
\end{lemma}

An analogous statement for algebras of uncountable cardinality can be obtained for a generalized quote, denoted $\chi,$ that has been defined in \cite{FontanellaGeoffroy2020}.

\section{Ordinals and Cardinals in \texorpdfstring{$\ZFepsilon$}{ZFepsilon}}

Over \tf{ZF}, there are multiple standard ways to define the notion of an ordinal. For example:
\begin{enumerate} \setlength \itemsep{3pt}
    \item An ordinal is a transitive set which is well-ordered by $\in$.
    \item An ordinal is a transitive set $a$ satisfying trichotomy. That is, for any $x, y, \in a$, precisely one of the following hold: $x \in y$, $y \in x$, or $x = y$.
    \item An ordinal is a transitive set of transitive sets.
\end{enumerate}

When working over the weaker theory $\ZFepsilon$, which does not have Extensionality, we should not necessarily expect the definitions to remain equivalent. We consider the following definitions. 


\begin{definition}\ 
    A set $a$ is said to be $\rlzin$\emph{-transitive} if $\forall x \rlzin a \forall y \rlzin x (y \rlzin a)$.

    A set $a$ is said to be an $\rlzin$\emph{-ordinal} if it is an $\rlzin$-transitive set of $\rlzin$-transitive sets.

    A set $a$ is said to be an \Rcomment{$\rlzin$\emph{-trichotomous ordinal}, or} $\rlzin$\emph{-TOD} if it is an $\rlzin$-transitive set which satisfies $\rlzin$\emph{-trichotomy}. That is, whenever $b, c \rlzin a$ then precisely one of the following hold: $b \rlzin c$, $c \rlzin b$, or $b = c$.
\end{definition}


We remark here that we have taken the definition of $\rlzin$-ordinal to be a transitive set of transitive sets which is not necessarily the standard definition. However, it turns out to be a much easier definition which a broad class of sets will satisfy. In particular, we will have that the $\rlzin$-ordinals form a proper class, whereas in general the class of $\rlzin$-TOD sets will be extensionally equal to a set.

\begin{proposition} \label{theorem:ETODiffEWOD}
    $(\ZFepsilon)$ The following are equivalent:
    \vspace{-5pt}
    \begin{enumerate} \setlength \itemsep{3pt}
        \item $a$ is an $\rlzin$-TOD.
        \item $a$ is an $\rlzin$-transitive set and $a$ is $\rlzin$-well founded, namely for every non-empty $X \subseteq_{\rlzin} a,$ $\exists z \rlzin X \, \forall x \rlzin X (z \rlzin x \lor z = x)$.
    \end{enumerate}
\end{proposition}

Before we give the proof, we first need two lemmas. The first gives us that in $\ZFepsilon$ if $x \rlzin y$ then $y \notrlzin x$. \Rcomment{This is usually proven using Foundation however, since this is not an assumed axiom of $\ZFepsilon$, we will use the axiom of $\rlzin$-induction instead.} The second tells us that any element of an $\rlzin$-TOD is itself an $\rlzin$-TOD. Note that the corresponding result is also trivially true for $\rlzin$-ordinals. 

\begin{lemma} \label{theorem:NoDecRlzinChains} \,
\vspace{-5pt}
    \begin{enumerate} \setlength \itemsep{3pt}
    \item $\ZFepsilon \vdash \forall x, y (x \rlzin y \rightarrow y \notrlzin x)$,
    \item $\ZFepsilon \vdash \forall x, y, z (x \rlzin y \land y \rlzin z \rightarrow z \notrlzin x)$.
    \end{enumerate}
\end{lemma}

\begin{proof}
    Both of these statements are proven using $\rlzin$-Induction. Since they are very similar we shall only prove the second of these. Let $\varphi(b) \equiv \forall x, y (x \rlzin b \land b \rlzin y \rightarrow y \notrlzin x)$ and suppose that $\forall b \rlzin a \varphi(b),$ we show $\varphi(a).$ Fix $b, c$ and suppose by contradiction that $b \rlzin a$, $a \rlzin c$ and $c \rlzin b$. Then, since $b \rlzin a$, we have $\varphi(b)$. So, since $c \rlzin b \land b \rlzin a$, we have $a \notrlzin c$, contradicting $a\rlzin c.$ Thus we must have $b \rlzin a \land a \rlzin c \rightarrow c \notrlzin a$, which is $\varphi(a)$. Thus we \Rcomment{have proven} $\ZFepsilon \vdash \forall \Rcomment{z} \varphi(\Rcomment{z})$.
\end{proof}

\begin{lemma}
    $(\ZFepsilon)$ Suppose that $a$ is an $\rlzin$-TOD and $b \rlzin a$, then $b$ is an $\rlzin$-TOD. 
\end{lemma}

\begin{proof}
    Suppose that $a$ is an $\rlzin$-TOD. Take $b \rlzin a$ and $d \rlzin c \rlzin b$, we need to show that $d \rlzin b$. First, suppose that $b = d$. Then we must have that $b \rlzin c \rlzin b$, which contradicts the first part of \Cref{theorem:NoDecRlzinChains}. On the other hand, suppose that $b \rlzin d$. Then, we would have that $b \rlzin d \rlzin c \rlzin b$ which contradicts the second part of \Cref{theorem:NoDecRlzinChains}. Hence, by the $\rlzin$-trichotomy of $a$, $d \rlzin b$.

    Finally, since $a$ is $\rlzin$-transitive, given $c, d \rlzin b$ we have $c, d \rlzin a$. Thus $c \rlzin d$ or $d \rlzin c$ or $c = d$.
\end{proof}

\begin{proof}[Proof of Proposition \ref{theorem:ETODiffEWOD}]
    We prove the first implication by $\rlzin$-Induction. So suppose that $a$ was an $\rlzin$-TOD and that every $b \rlzin a$ satisfied (2). Now suppose for a contradiction that for some $X \subseteq_{\rlzin} a$, $\forall z \rlzin X \exists x \rlzin X \neg (z \rlzin x \lor z = x)$. Since $a$ is $\rlzin$ trichotomous, it must then be the case that $\forall z \rlzin X \exists x \rlzin X (x \rlzin z)$. Given $b \rlzin X$ let $Y = \{ y \rlzin X \divline y \rlzin b \}$, which is non-empty by assumption. Since $Y \subseteq_{\rlzin} b$, which satisfies (2), we can fix some $z \rlzin Y$ such that $\forall x \rlzin Y (z \rlzin x \lor z = x)$. Now, consider $x \rlzin X$. If $x \rlzin z$ then, since $z \rlzin b$ and $b$ is $\rlzin$-transitive, $x \rlzin Y$. But this means that $z \rlzin x \lor z = x$, which gives a contradiction. Hence, since $x, z \rlzin a$, which is $\rlzin$-trichotomous, we must have $z \rlzin x \lor z = x$, contradicting the assumption that $X$ did not have an $\rlzin$-minimal element.

    For the reverse implication, suppose that $a$ satisfied (2) and take $x \neq y$ in $a$. Consider $X = \{ z \rlzin a \divline z = x \lor z = y \lor (z \rlzin x \land z \notrlzin y) \lor (z \rlzin y \land z \notrlzin x) \}$. This has an $\rlzin$-minimal element, say $z$. Now, if $z \rlzin x$ but $z \notrlzin y$ then, by minimality, we must have $z = y$. Hence $y \rlzin x$. Similarly, if $z \rlzin y$ but $z \notrlzin x$, then $z = x \rlzin y$. Thus $x \rlzin y \lor y \rlzin x$, so $a$ is $\rlzin$-trichotomous. 
\end{proof}

\begin{corollary}
    $(\ZFepsilon)$ Suppose that $a$ is an $\rlzin$-TOD. Then $a$ is an $\rlzin$-ordinal.
\end{corollary}

While every $\rlzin$-TOD is also an $\rlzin$-ordinal the reverse need not hold. For example, suppose that $\fullname{2}$ had size greater than $2$. Then there will be two ordinals $a, b \rlzin \fullname{2}$ such that $a \notrlzin b$ and $b \notrlzin a$. But then $(a \cup \{a\}) \cup (b \cup \{b\})$ is an $\rlzin$-ordinal which is not $\rlzin$-trichotomous. In the same way, it is possible to have multiple $\rlzin$-\Rcomment{ordinals} which externally have the same size\Rcomment{. For example, in \cite{FontanellaGeoffroy2020} the authors introduce the ordinals $\hat{n}$ (which we will define in the next section) and show that $\hat{n}$ is an $\rlzin$-TOD which is extensionally equal to $n$. From this it follows that $\hat{4}$ always has size $4$, while it is also possible for $\fullname{2}$ to have size $4$, such as in \cite{Krivine2018}.} \\

A final, key property, of $\rlzin$-TODs is that \Rcomment{extensional} equality is equivalent to Leibniz equality.

\begin{proposition}(\cite[Proposition 3.10]{FontanellaGeoffroy2020})
    $(\ZFepsilon)$ Suppose that $a$ is an $\rlzin$-TOD. Then $\forall x \rlzin a \, \forall y \rlzin a (x \simeq y \rightarrow x = y)$.
\end{proposition}

\begin{proof}
    Let $a$ be an $\rlzin$-TOD and take $x, y \rlzin a$ with $x \neq y$. Then either $x \rlzin y$ or $y \rlzin x$. But, if $x \rlzin y$ then $x \in y$, from which it follows that $x \not\simeq y$. Similarly, if $y \rlzin x$ then $y \not\simeq x$.
\end{proof}

We can also define the concept of an $\rlzin$-cardinal over $\ZFepsilon$, which we do next.

\begin{definition}
An ordinal $a$ is said to be an $\rlzin$\emph{-cardinal} if for every $b \rlzin a$ there is no extensional function which is an $\rlzin$-surjection of $b$ onto $a$.
\end{definition}

It may appear unnatural to define $\rlzin$-cardinals in terms of extensional functions rather than $\rlzin$-functions. However, this is done because we want to have that if $a$ is an $\rlzin$-cardinal in a model of $\ZFepsilon$ then it is an $\in$-cardinal in the extensional part. This then allows us to obtain results in the resulting $\tf{ZF}$ model.

\begin{proposition} $(\ZFepsilon)$
\begin{enumerate}
    \item If $a$ is an $\rlzin$-transitive set, then it is an $\in$-transitive set;
    \item If $a$ is an $\rlzin$-ordinal, then it is an $\in$-ordinal;
    \item If $a$ is an $\rlzin$-cardinal, then it is an $\in$-cardinal.
\end{enumerate}
\end{proposition}

\begin{proof} \,
\begin{enumerate}
    \item Take $c \in b \in a$. Then there exists some $x \rlzin a$ such that $x \simeq b$ and there exists some $y \rlzin x$ such that $y \simeq c$. Since $a$ is assumed to be $\rlzin$-transitive, $y \rlzin a$. Therefore $c \in a$ by definition of $\in$.
    \item Suppose that $a$ is an $\rlzin$-ordinal. We have already shown that $a$ is $\in$-transitive, so it suffices to prove that every $b \in a$ is $\in$-transitive. So let $d \in c \in b \in a$. Then we can find $z \rlzin y \rlzin x \rlzin a$ such that $x \simeq b$, $y \simeq c$ and $z \simeq d$. Since $a$ is an $\rlzin$-ordinal, $z \rlzin x$ and therefore $d \in x$. Finally, $d \in x$ and $x \simeq b$ gives us $d \in b$, as required.
    \item Let $a$ be an $\rlzin$-ordinal. We shall show that if $a$ is not an $\in$-cardinal then $a$ is not an $\rlzin$-cardinal. To this end, suppose that $b \in a$ and $f \colon b \rightarrow a$ was an $\in$-surjection. Fix $z \rlzin a$ such that $b \simeq z$ and define a binary relation $g$ by 
    \[
    g(x, y) \Longleftrightarrow \exists c \in a \, \big( f(x, c) \land c \simeq y \big).
    \]
    We shall show that $g$ defines an extensional function which is an $\rlzin$-surjection from $z$ onto $a$. There are three things to prove: that $g$ is $\rlzin$-total on $a$; $g$ satisfies the definition of being an extensional function; and $g$ is an $\rlzin$-surjection.

    For the first claim, take $x \rlzin z$. Then $x \in z$ and therefore $x \in b$. Since $f$ is $\in$-total, we can fix $c \in a$ such that $f(x, c)$ and then fix $y \rlzin a$ such that $y \simeq c$. Hence we have $g(x,y)$.

    For the second claim, fix $x, x' \rlzin b$, $y, y' \rlzin a$ and suppose that $g(x, y)$ holds, $g(x', y')$ holds and $x \simeq x'$. Fix $c, c' \in a$ such that $f(x, c)$ holds, $f(x',c')$ holds, $c \simeq y$ and $c' \simeq y'$. Since $f$ is an $\in$-function, $x, x' \in b$ and $x \simeq x'$, we must have that $c \simeq c'$. Therefore $y \simeq c \simeq c' \simeq y'$.

    For the final claim, fix $y \rlzin a$. Then $y \in a$ so, since $f$ is an $\in$-surjection, we can fix $t \in b$ such that $f(t, y)$ holds. Since $t \in b$ and $b \simeq z$, $t \in z$. Therefore, there exists some $d \rlzin z$ such that $t \simeq d$. Since $g$ is $\rlzin$-total, we can fix $y' \rlzin a$ such that $g(d, y')$ holds and then, by definition, $c \in a$ such that $f(d, c)$ and $c \simeq y'$. Finally, since $f$ is an $\in$-function and $t \simeq d$, $y \simeq c \simeq y'$, from this we can conclude that $g(d, y)$ holds.
\end{enumerate}
\end{proof}

In our arguments, it will often be beneficial to improve the claim that there are no surjections to one asserting that every function is bounded. This will be particularly useful because it allows us to work with the notion of $\rlzin$-functions rather than extensional ones. On the other hand, we note here that the non-extensional axiom of choice is vital to the proof.

\begin{proposition} \label{theorem:NEACGivesCardinals}
    $(\ZFepsilon + \tf{NEAC})$ Suppose that $a$ is an $\rlzin$-ordinal and for all $b \rlzin a$, every $\rlzin$-function $f \colon b \rightarrow a$ is bounded, that is to say there is some $c \rlzin a$ such that $f \pointwise b \subseteq_{\rlzin} c$. Then $a$ is an $\in$-cardinal.
\end{proposition}

\begin{proof}
    We shall show that for every $b \in a$, every $\in$-function $f \colon b \rightarrow a$ is bounded, and therefore not an $\in$-surjection.

    To do this, take $b \in a$ and $b' \rlzin a$ with $b \simeq b'$. Now, suppose that $f \colon b \rightarrow a$ were an $\in$-function. Define a relation $r$ on $b' \times a$ by 
        \[
        \op(x, y) \rlzin r \Longleftrightarrow f(x) \simeq y.
        \]
    
    Then, by \tf{NEAC}, we can find an $\rlzin$-function $g \subseteq_{\rlzin} r$ with $g \colon b' \rightarrow a$. Since $g$ is an $\rlzin$-function it is bounded, say by $c \rlzin a$. Now, if $x \in b$, there is some $x' \rlzin b'$ with $x \simeq x'$. Then $f(x) \simeq f(x') \simeq g(x') \rlzin c$. Thus $f(x) \in c$, so $f \pointwise b \subseteq c$ and $f$ is bounded.
\end{proof}

\section{Names for Ordinals}

Now, given an ordinal $\alpha$ in the ground model we have two ways to associate it to an ordinal in a realizability model: $\fullname{\alpha} \Rcomment{ = \{ (\fullname{\beta}, \pi) \divline \beta \in \alpha, \pi \in \Pi \}}$ and $\hat{\alpha}$ (the former is presented in \cite{Matthews2023}, the latter was introduced in \cite{FontanellaGeoffroy2020}), which we discuss in this section. The key difference between them is that $\hat{\alpha}$ will be an $\rlzin$-TOD whereas $\fullname{\alpha}$ will only be if the realizability \Rcomment{model} is trivial. On the other hand, we can only define $\hat{\alpha}$ for $\alpha < | \Lambda|$ and to prove that it is an $\rlzin$-TOD will require the existence of a special realizer (see \cite{FontanellaGeoffroy2020} for details). Instead, $\fullname{\alpha}$ is defined for every ordinal $\alpha$ \Rcomment{in} any realizability model without any need for additional instructions.

\begin{definition}
    Suppose that $\langle \nu_\alpha \divline \alpha < | \Lambda | \rangle$ is an enumeration of all terms. For $\alpha \leq | \Lambda |$ define $\hat{\alpha}$ recursively as
    \[
    \hat{\alpha}  = \{ (\hat{\beta}, \nu_\beta \stackapp \pi ) \divline \beta \in \alpha\Rcomment{, \pi \in \Pi} \}.
    \]
\end{definition}

We state here some of the standard properties one can prove about these names.

\begin{proposition} \,
    \begin{enumerate}
        \item For all $\beta < \alpha$, $\identity \Vdash \fullname{\beta} \rlzin \fullname{\alpha}$ and $\lamAbstOne \lapp{\app{u}{\rlzXSubsetX}}{\rlzXSubsetX} \Vdash \Rcomment{\fullname{\beta} \subseteq \fullname{\alpha}}$,
        \item For all $\beta < \alpha \leq | \Lambda |$, $\lamAbstOne \app{\lamtermOne}{\nu_\beta} \Vdash \hat{\beta} \rlzin \hat{\alpha}$ and $\lamAbstOne \lapp{\app{u}{\rlzXSubsetX}}{\rlzXSubsetX} \Vdash \Rcomment{\hat{\beta} \subseteq \hat{\alpha}}$.
    \end{enumerate}
\end{proposition}

We can also define a variant of bounded quantification which works for restrictions to hat names. 

\begin{definition}
    For $\alpha \leq | \Lambda |$ and $\varphi(x)$ a formula, let
    \[
    \falsity{\forall x^{\hat{\alpha}} \varphi(x)} = \bigcup_{\beta \in \alpha} \{ \nu_\beta \stackapp \pi \divline \pi \in \falsity{\varphi(\hat{\beta})}.
    \]
\end{definition}

As before, it is easy to see that this exhibits the same behaviour as bounded quantification. \Lcomment{The notation for the bounded quantifier $\forall x^{\hat{\alpha}}$ may be confused with the bounded quantifier $\forall x^{\cjgimel(A)}$ defined in \Cref{definition:BoundedUniversals}, but it will be clear from the context which of the two we mean.}

\begin{proposition}
    For every $\alpha \leq | \Lambda |$,
    \begin{enumerate}
        \item $\lamAbstOne \lamAbstTwo \lamAbstThree \rapp{v}{\app{u}{w}} \Vdash \forall x^{\hat{\alpha}} \varphi(x) \rightarrow \forall x (\neg \varphi(x) \rightarrow x \notrlzin \hat{\alpha})$,
        \item $\lamAbstOne \app{\cc}{u} \Vdash \forall x (\neg \varphi(x) \rightarrow x \notrlzin \hat{\alpha}) \rightarrow \forall x^{\hat{\alpha}} \varphi(x)$.
    \end{enumerate}
\end{proposition}

We next show that both $\fullname{\alpha}$ and $\hat{\alpha}$ correspond to ordinals in the realizability model. We shall observe here that so far there was no need for additional special instructions; these will be needed to show that $\hat{\alpha}$ is an $\rlzin$-TOD.

\begin{proposition} \label{DalethAlphaOrdinal}
If $\delta$ is an ordinal in \tf{V} then it is realized that $\fullname{\delta}$ is an $\rlzin$-ordinal in $\rlzmodel$ by a realizer that does not depend on $\delta$.
\end{proposition}

\begin{proof}
Let $\delta$ be an ordinal in \tf{V}. We shall show that $\fullname{\delta}$ is an $\rlzin$-transitive set of $\rlzin$-transitive sets. To do this it suffices to show that
\[
\identity \Vdash \forall x^{\fullname{\delta}} \forall y (y \notrlzin \fullname{\delta} \rightarrow y \notrlzin x)
\]
and
\[
\lamAbstOne \lamAbstTwo \app{\lamtermTwo}{\lamtermOne} \Vdash \forall x^{\fullname{\delta}} \forall y \forall z (z \notrlzin x \rightarrow ( z \rlzin y \rightarrow y \notrlzin x)).
\]
For the first claim, fix $\beta \in \delta$, $c \in \rlzstr$, $t \Vdash c \notrlzin \fullname{\delta}$ and $\pi \in \falsity{c \notrlzin \fullname{\beta}}$. Now $\falsity{c \notrlzin \fullname{\beta}} = \{ \sigma \divline (c, \sigma) \in \fullname{\beta} \}$. So, since this set is non-empty, it must be the case that $\falsity{c \notrlzin \fullname{\beta}} = \Pi$ and $c = \fullname{\gamma}$ for some $\gamma \in \beta$. Therefore, $\falsity{c \notrlzin \fullname{\delta}} = \falsity{\fullname{\gamma} \notrlzin \fullname{\delta}} = \Pi$ and therefore $t \star \pi \in \Perp$, from which the result follows.

For the second claim, fix $\beta < \alpha < \delta$, $c \in \rlzstr$, $t \Vdash c \notrlzin \fullname{\alpha}$, $s \Vdash c \rlzin \fullname{\beta}$ and $\pi \in \Pi$. We first show that $t \Vdash c \notrlzin \fullname{\beta}$.

If $c \neq \fullname{\gamma}$ for any $\gamma \in \beta$ then $\falsity{c \notrlzin \fullname{\beta}} = \emptyset$ and the result follows so suppose that $c = \fullname{\gamma}$ for some $\gamma \in \beta$. Then, since $\beta \subseteq \alpha$, $\falsity{c \notrlzin \fullname{\beta}} = \falsity{c \notrlzin \fullname{\alpha}} = \Pi$. Therefore $t \star \sigma \in \Perp$ for any $\sigma \in \Pi$. Hence $t \star \sigma \in \Perp$ for every $\sigma \in \falsity{c \notrlzin \fullname{\beta}}$.

From this it follows that $\lamAbstOne \lamAbstTwo \app{\lamtermTwo}{\lamtermOne} \star t \stackapp s \stackapp \pi \succ s \star t \stackapp \pi \in \Perp$.
\end{proof}

\begin{proposition} \label{HatAlphaOrdinal}
    If $\delta \leq | \Lambda|$ then it is realized that $\hat{\delta}$ is an $\rlzin$-ordinal in $\rlzmodel$ by a realizer that does not depend on $\delta$.
\end{proposition}

\begin{proof}
    It will suffice to prove that
    \[
    \lamAbstOne \lamAbstThree \app{\lamtermThree} \Vdash \forall x^{\Rcomment{\hat{\delta}}} \forall y (y \notrlzin \hat{\delta} \rightarrow y \notrlzin x) 
    \]
    and
    \[
    \lamAbstOne \lamAbstTwo \lamAbstThree \app{\lamtermThree}{\lamtermTwo} \Vdash \forall x^{\hat{\delta}} \forall y \forall z ( z \notrlzin x \rightarrow ( z \rlzin y \rightarrow y \notrlzin x)).
    \]
    For the first claim, $\falsity{ \forall x^{\Rcomment{\hat{\delta}}} \forall y (y \notrlzin \hat{\delta} \rightarrow y \notrlzin x)} = \bigcup_{\alpha < \beta < \delta} \{ \nu_\beta \stackapp t \stackapp \nu_\alpha \stackapp \pi \divline t \Vdash \hat{\alpha} \notrlzin \hat{\delta}, \, \pi \in \Pi \}$. So, fix $\alpha < \beta < \delta$, $\pi \in \Pi$ and suppose $t \Vdash \hat{\alpha} \notrlzin \hat{\delta}$. Then 
    \[
    \lamAbstOne \lamAbstThree \lamtermThree \star \nu_\beta \stackapp t \stackapp \nu_\alpha \stackapp \pi \succ t \star \nu_\alpha \stackapp \pi \in \Perp.
    \]
    For the second claim, observe that element of $\falsity{\forall x^{\hat{\delta}} \forall y \forall z ( z \notrlzin x \rightarrow ( z \rlzin y \rightarrow y \notrlzin x))}$ are of the form $\nu_\gamma \stackapp t \stackapp s \stackapp \nu_\beta \stackapp \pi$ where $\beta < \gamma < \delta$, $a \in \rlzstr$, $t \Vdash a \notrlzin \hat{\gamma}$, $s \Vdash a \rlzin \hat{\beta}$ and $(\hat{\beta}, \nu_\beta \stackapp \pi) \in \hat{\gamma}$. We shall first see that $t \Vdash a \notrlzin \hat{\beta}$.

    For this, if $a \neq \hat{\alpha}$ for any $\alpha < \beta$ then $\falsity{c \notrlzin \hat{\beta}} = \emptyset$ and the result is immediate. So suppose that $a = \hat{\alpha}$ for some $\alpha < \beta$. Then $\falsity{c \notrlzin \hat{\beta}} = \{ \nu_\alpha \stackapp \sigma \divline \sigma \in \Pi \} = \falsity{c \notrlzin \hat{\gamma}}$. Thus $s \star t \stackapp \pi' \in \Perp$ for all $\pi' \in \Pi$ and hence
    \[
    \lamAbstOne \lamAbstTwo \lamAbstThree \app{\lamtermThree}{\lamtermTwo} \star \nu_\gamma \stackapp t \stackapp s \stackapp \nu_\beta \stackapp \pi \succ s \star t \stackapp \nu_\beta \stackapp \pi \in \Perp. 
    \]
\end{proof}

Since $\ZFepsilon$ proves that any $\rlzin$-ordinal is an $\in$-ordinal and any $\in$-ordinal satisfies trichotomy, and the above realizers were uniform, the following corollary is immediate.

\begin{corollary}
    There exist fixed realizers $u_0$ and $u_1$ such that for any $\delta \in \tf{Ord}$ and $\mu \in | \Lambda |$,
    \[
    u_0 \Vdash \forall x^{\fullname{\delta}} \forall y^{\fullname{\delta}} ( x \in y \lor x \simeq y \lor y \in x)
    \]
    and    
    \[
    u_1 \Vdash \forall x^{\hat{\mu}} \forall y^{\hat{\mu}} ( x \in y \lor x \simeq y \lor y \in x)
    \]
\end{corollary}

We end with the observation that if $\alpha$ is a limit ordinal then $\fullname{\alpha}$ is \Rcomment{also} realized to be \Rcomment{an $\rlzin$-limit ordinal. Recall that $\alpha$ is a limit ordinal if for every $\beta \in \alpha$, $\beta + 1 \in \alpha$. However, since there are many non-extensionally equal names for $\beta + 1$ we need to take a slightly weaker formulation of being a $\rlzin$-limit ordinal.}

\begin{definition}[$\ZFepsilon$]
    We say that an $\rlzin$-ordinal $a$ is \Rcomment{an $\rlzin$-limit ordinal} if $\forall x \rlzin a \exists y \rlzin a (x \rlzin y)$.
\end{definition}

\begin{proposition}[$\ZFepsilon$]
    If $a$ is \Rcomment{an $\rlzin$-limit ordinal} then $a$ is \Rcomment{an $\in$-limit ordinal}.
\end{proposition}

\begin{proof}
    Suppose that $a$ is an $\rlzin$-ordinal (and hence an $\in$-ordinal) and that $\forall x \rlzin a \exists y \rlzin a (x \rlzin y)$. Fix $x \in a$ and $x' \rlzin a$ such that $x \simeq x'$. Then we can find some $y \rlzin a$ such that $x' \rlzin y$. But from this it follows that $x \in y$. Thus we have that $\forall x \in a \exists y \in a (x \in y)$.

    \Rcomment{To see that this implies that $a$ is an $\in$-limit ordinal}, suppose that $\alpha \in a$ and take $\beta \in a$ such that $\alpha \in \beta$. Then, since the $\in$-ordinals are linearly ordered by $\in$, $\alpha + 1 \leq \beta$. So, since $\beta \in a$ and $a$ is an ordinal we must have that $\alpha + 1 \in a$.
\end{proof}

\begin{proposition} \label{theorem:DalethPreservesLimits}
    If $\alpha$ is a limit ordinal then $\rlzmodel \Vdash ``\fullname{\alpha} \text{ is an } \Rcomment{\rlzin\text{-limit ordinal''}}$. Hence $\fullname{\alpha}$ is realized to be \Rcomment{an $\in$-limit ordinal}.
\end{proposition}

\begin{proof}
    Let $\alpha$ be a limit ordinal. It suffices to prove that $\identity \Vdash \forall x^{\fullname{\alpha}} (\forall y^{\fullname{\alpha}} (x \notrlzin y) \rightarrow \perp).$ To do this, fix $\beta \in \alpha$, $\pi \in \Pi$ and suppose that $t \Vdash \forall y^{\fullname{\alpha}} (\fullname{\beta} \notrlzin y)$. Since $\alpha$ is a limit ordinal, $\beta + 1 \in \alpha$ and thus $t \Vdash \fullname{\beta} \notrlzin \fullname{\beta + 1}$. But, $\falsity{\fullname{\beta} \notrlzin \fullname{\beta + 1}} = \Pi$ and thus $\identity \star t \stackapp \pi \succ t \stackapp \pi \in \Perp$.
\end{proof}

Finally,  in order to obtain that $\hat{\alpha}$ is an $\rlzin$-TOD we need to add a special instruction, $\chi$, which allows one to compare the terms in $\Lambda$ by their indices. This is done in section 3 of \cite{FontanellaGeoffroy2020}.

\begin{definition} \label{definition:PropertyChi}
    Let $\chi \in \Lambda$ be a special instruction defined by the following rule: we extend the order $\prec$ to be the smallest pre-order on $\Lambda \star \Pi$ such that for any $\alpha, \beta \in | \Lambda |$, $t, s, r \in \Lambda$ and $\pi \in \Pi$,
    \[
    \chi \star \nu_\alpha \stackapp \nu_\beta \stackapp t \stackapp s \stackapp r \stackapp \pi \succ \begin{cases}
        t \star \pi & \text{if }\ \alpha < \beta, \\
        s \star \pi & \text{if }\ \alpha = \beta, \\
        r \star \pi & \text{if }\ \beta < \alpha.
    \end{cases}
    \]
\end{definition}

\begin{theorem}(\cite[Proposition 3.11]{FontanellaGeoffroy2020})
    Suppose that $\chi \in \Lambda$ and $\prec$ has been expanded to satisfy \Cref{definition:PropertyChi}. Then there exists a realizer $\theta \in \mathcal{R}$ such that for every $\alpha \leq | \Lambda |$,
    \[
    \theta \Vdash ``\hat{\alpha} \text{ is \Rcomment{an} } \rlzin\text{-Trichotomous ordinal''}.
    \]
\end{theorem}

\section{Comparing Names for Ordinals}

Having defined the names $\fullname{\alpha}$ and $\hat{\alpha}$, we now wish to see which ordinals these names correspond to in a realizability model. It is easy to see that $\fullname{0} = \hat{0} = \emptyset$, which will be a name for the empty set in $\rlzmodel$. In this section we investigate three key things: successor ordinals, $\omega$ (the first limit ordinal), the class of all ordinals.

\subsection{Computing Successors}

It has \Rcomment{previously} been proven in \cite{FontanellaGeoffroy2020} that $\widehat{n + 1}$ extensionally corresponds to $n + 1$ and $\hat{\omega}$ corresponds to $\omega$. \Rcomment{Using an alternative proof, we shall reprove this and also} show that the same holds for the recursive names. For this, recall that $\rlzXSubsetX$ was a fixed realizer such that $\rlzXSubsetX \Vdash \forall x (x \subseteq x)$.

\begin{lemma} \label{theorem:DalethnLargestElement}
    For any ordinal $\alpha$, 
    \[
    \lamAbstOne \lapp{\app{\lamtermOne}{\rlzXSubsetX}}{\rlzXSubsetX} \Vdash \forall x^{\fullname{\alpha + 1}} ( x \subseteq \fullname{\alpha}).
    \]
\end{lemma}

\begin{proof}
    Fix $\beta \in \alpha + 1$ and $t \stackapp \pi \in \falsity{\fullname{\beta} \subseteq \fullname{\alpha}}$. This means that for some $\gamma \in \beta$ we have $t \Vdash \fullname{\gamma} \not\in \fullname{\alpha}$ while $(\fullname{\gamma}, \pi) \in \fullname{\beta}$. Since $\gamma < \beta \leq \alpha$, $(\fullname{\gamma}, \pi) \in \fullname{\alpha}$ and thus $\rlzXSubsetX \stackapp \rlzXSubsetX \stackapp \pi \in \falsity{\fullname{\gamma} \not\in \fullname{\alpha}}$. Hence $\lamAbstOne \lapp{\app{\lamtermOne}{\rlzXSubsetX}}{\rlzXSubsetX} \star t \stackapp \pi \succ t \star \rlzXSubsetX \stackapp \rlzXSubsetX \stackapp \pi \in \Perp$.
\end{proof}

From this it follows that if $\alpha$ is a successor ordinal then so is $\fullname{\alpha}$ \Rcomment{in the extensional sense}. \Rcomment{We remark here that we shall make free use of the symbols $+$, $<$ and $\leq$} \Rcomment{when discussing extensional ordinals. So, by $\beta \leq \alpha$ we mean $\beta \in \alpha$ or $\beta \simeq \alpha$ and by $\alpha + 1$ we mean $\alpha \cup \{\alpha\}$.}

\begin{proposition} \label{theorem:RecursiveSuccessors}
    If $\alpha = \beta + 1$ then $\rlzmodel \Vdash ``\fullname{\alpha} \text{ is the } \Rcomment{\in}\text{-successor of } \fullname{\beta}\text{''}$, that is \hbox{$\rlzmodel \Vdash \fullname{\alpha} \simeq \fullname{\beta} + 1$.}
\end{proposition}

\begin{proof}
   Since $\beta \in \alpha$, $\identity \Vdash \fullname{\beta} \rlzin \fullname{\alpha}$. Thus, using \Cref{theorem:DalethnLargestElement}, 
   \[
   \rlzmodel \Vdash \fullname{\beta} \rlzin \fullname{\alpha} \land \forall x^{\fullname{\alpha}} (x \subseteq \fullname{\beta}).
   \]
   We now argue within a realizability model $\rlzmodel = (\rlzstr, \notrlzin, \not\in, \subseteq)$. We know that $\fullname{\beta}$ and $\fullname{\alpha}$ are $\in$-ordinals with $\fullname{\beta} < \fullname{\alpha}.$ Thus $\fullname{\beta} + 1 \leq \fullname{\alpha}$ (here ``$a + 1$'' is some set which is extensionally equal to the successor of the $\in$-ordinal $a$). On the other hand, \Rcomment{using \Cref{theorem:DalethnLargestElement}, \Cref{theorem:BoundedUniversals} and the fact that $x \rlzin a \rightarrow x \in a$}, we also have that for all $x \in \fullname{\alpha}$, $x$ is an $\in$-ordinal with $x \leq \fullname{\beta}$. From this it follows that $\fullname{\alpha} \leq \fullname{\beta} + 1$ and thus the two are indeed extensionally equal.
\end{proof}

\begin{corollary}
    For all $n \in \omega$, $\rlzmodel \Vdash \fullname{n} \text{ is the } n\textsuperscript{th}\text{-successor of } \emptyset$.
\end{corollary}

In order to prove that $\widehat{n + 1}$ is extensionally equal to the successor of $\hat{n}$, in \cite{FontanellaGeoffroy2020} the authors need to add another additional instruction $\xi$. This instruction was used to ensure that successor ordinals do not ``collapse'' to their predecessor. However, by a very similar argument to the above, we can show that this additional instruction can be removed. In particular, in \emph{any} realizability model we will have that $\widehat{\alpha + 1}$ names the successor of $\hat{\alpha}$.\footnote{Note that this argument will also not require the special instruction $\chi$. This was only needed to ensure that $\hat{\alpha}$ was an $\rlzin$-TOD.}

For this, we give a mild variant of the theory of lifting functions. This will give us a method to lift any ground model function $f \colon \kappa \rightarrow \kappa$, where $\kappa = | \Lambda|$, to an $\rlzin$-function $\hat{f} \colon \hat{\kappa} \rightarrow \hat{\kappa}$ in the realizability model.

\begin{definition}
    Let $\kappa = | \Lambda |$ and $f \colon \kappa \rightarrow \kappa$. We define the \emph{ordered lift} of $f$, $\hat{f}$ as
    \[
    \hat{f} \coloneqq \{ ( \op(\hat{\alpha}, \widehat{f(\alpha)}), \nu_\alpha \stackapp \pi) \divline \alpha \in \kappa, \pi \in \Pi \}.
    \]
\end{definition}

As before, we can also define the value $\falsity{\varphi(x, \hat{f}(x))}$ analogously to the general case.

\begin{definition}
    Let $f \colon \kappa \rightarrow \kappa$ be a function in the ground model. Then, for any formula $\varphi(u, v, \overrightarrow{w})$ and $\overrightarrow{w} \subset \rlzstr$,
    \[
    \falsity{\varphi(x, \hat{f}(x), \overrightarrow{w})} = \begin{cases}
        \{ \nu_\alpha \stackapp \pi \divline \pi \in \falsity{\varphi(\hat{\alpha}, \widehat{f(\alpha)}, \overrightarrow{w})} \Rcomment{\}} & \text{ if } x = \alpha \text{ for some } \alpha \in \kappa, \\
        \emptyset & \text{ otherwise.}
    \end{cases}
    \]
\end{definition}

As for the traditional lift, it is straightforward to show that $\hat{f}$ is a lift of the function $f$.

\begin{proposition}
    If $f \colon \kappa \rightarrow \kappa$ is any function \Rcomment{in} the ground model then, for every formula $\varphi$,
    \[
    \rlzmodel \Vdash \forall x^{\hat{\kappa}} \big( \varphi(x, \hat{f}(x), \overrightarrow{w}) \longleftrightarrow \forall y ( \op(x, y) \rlzin \hat{f} \rightarrow \varphi(x, y, \overrightarrow{w})) \big).
    \]
\end{proposition}

\begin{proposition} \label{theorem:HatLiftProperties}
    If $f \colon \kappa \rightarrow \kappa$ is any function is the ground model then, in $\rlzmodel$, $\hat{f}$ is an $\rlzin$-function with domain 
    \Rcomment{$\hat{\kappa}$}
    and range $\hat{\kappa}$. Moreover, for any $\alpha \in \kappa$,
    \[
    \rlzmodel \Vdash \hat{f}(\hat{\alpha}) = \widehat{f(\alpha)}.
    \]
\end{proposition}


We now consider the successor function $\ff{succ} \colon \kappa \rightarrow \kappa$, $\ff{succ}(\alpha) = \alpha + 1.$ 

\begin{lemma} \label{theorem:HatLiftIsSuccessor}
    Let $\rlzXSubsetX$ be the fixed realizer such that $\rlzXSubsetX \Vdash \forall x (x \subseteq x)$. Then 
    \[
    \lamAbstTwo \lamAbstThree \rapp{\lamtermThree}{\lamAbstOne \lapp{\app{u}{\rlzXSubsetX}}{\rlzXSubsetX}} \Vdash \forall x^{\hat{\kappa}} \forall y (( y \subseteq x \rightarrow \perp) \rightarrow y \notrlzin \widehat{\ff{succ}}(x)).
    \]
    Hence, $\rlzmodel \Vdash \forall x^{\hat{\kappa}} ( \widehat{\ff{succ}}(x) \simeq x + 1)$.
\end{lemma}

\begin{proof}
    First observe that any element of $\falsity{\forall x^{\hat{\kappa}} \forall y (( y \subseteq x \rightarrow \perp) \rightarrow y \notrlzin \widehat{\ff{succ}}(x))}$ is of the form $\nu_\alpha \stackapp t \stackapp \nu_\beta \stackapp \pi$ where $\beta < \alpha < \kappa$, $t \Vdash \hat{\beta} \subseteq \hat{\alpha} \rightarrow \perp$ and $(\nu_\beta, \pi) \in \widehat{\ff{succ}}(\hat{\alpha}) = \widehat{\alpha + 1}$. But, since $\beta \leq \alpha$, $\lamAbstOne \lapp{\app{u}{\rlzXSubsetX}}{\rlzXSubsetX} \Vdash \hat{\beta} \subseteq \hat{\alpha}$. Thus, $t \star (\lamAbstOne \lapp{\app{u}{\rlzXSubsetX}}{\rlzXSubsetX}) \stackapp \sigma \in \Perp$ for all $\sigma \in \Pi$, from which the result follows.
\end{proof}

\begin{corollary}
    For all $n \in \omega$, $\rlzmodel \Vdash \fullname{n} \simeq \hat{n}$.
\end{corollary}

\subsection{The ordinals of \texorpdfstring{$\rlzmodel$}{N}}

\Lcomment{In this section we show that the class of all ordinals of $\rlzmodel$ corresponds to the class $\fullname{\tf{Ord}} = \{ (\fullname{\alpha}, \pi) \divline \alpha \in \tf{Ord}, \pi \in \Pi \}$. Note that we can expand the language of realizability structure with predicates for arbitrary classes $C$ by extending the definition of the truth and falsity values in the obvious way for formulas of the form $x\rlzin C,$ $x\in C$ and $x \subseteq C$.}
We can define \Lcomment{a} bounded operator $\forall x^{\fullname{\tf{Ord}}}$ by
\[
\falsity{\forall x^{\fullname{\tf{Ord}}} \varphi(x)} = \bigcup_{\alpha \in \tf{Ord}} \falsity{\varphi(\fullname{\alpha}}
\]
and it is easy to see that this has the desired meaning.

\begin{lemma}
    $\rlzmodel \Vdash \forall x^{\fullname{\tf{Ord}}} ( \tf{Ord}_{\rlzin}(x) )$, where $\tf{Ord}_{\rlzin}(x)$ is the assertion that $x$ is an $\rlzin$-transitive set of $\rlzin$-transitive sets.
\end{lemma}

\begin{proof}
    This just follows from the fact that there is a realizer for the statement $\fullname{\alpha} \text{ is a } \rlzin\text{-transitive set of } \rlzin\text{-transitive sets}$ which does not depend on $\alpha$.
\end{proof}

We show that $\fullname{\tf{Ord}}$ is a proper class. 

\begin{lemma}\label{lemma: reish ord is proper class}
    For any set $a \in \rlzstr$, there exists an ordinal $\beta_a$ such that
    \[
    \lamAbstOne \lamAbstTwo \lapp{\lapp{\app{\rlzSubsetTrans}{\lamtermOne}}{\lamtermTwo}}{\rlzXNotInX} \Vdash \fullname{\beta_a} \not\in a,
    \]
    where $\rlzXNotInX, \rlzSubsetTrans \in \mathcal{R}$ are the fixed realizers such that $\rlzXNotInX \Vdash \forall x (x \not\in x)$ and $\rlzSubsetTrans \Vdash \forall x \forall y \forall z (x \subseteq y \rightarrow (y \subseteq z \rightarrow x \subseteq z)).$
\end{lemma}

\begin{proof}
    Given $c \in \ff{dom}(a)$ and $s \in \Lambda$, define
    \[
    \delta_{c, s} \coloneqq 
    \begin{cases}
        \ff{min}\{ \delta \divline s \Vdash c \subseteq \fullname{\delta} \} & \text{if this set is nonempty} \\
        \emptyset & \text{otherwise}.
    \end{cases}
    \]
    Next, let $\beta_a \coloneqq \ff{sup}\{ \delta_{c, s} \divline c \in \ff{dom}(a), \, s \in \Lambda \} + 1$.

    Now, suppose that $t \stackapp s \stackapp \pi \in \falsity{\fullname{\beta_a} \not\in a}$. Then we can fix $c$ such that $(c, \pi) \in a$, $t \Vdash \fullname{\beta_a} \subseteq c$ and $s \Vdash c \subseteq \fullname{\beta_a}$. By definition, this means that there exists some $\delta < \beta_a$ such that $s \Vdash c \subseteq \fullname{\delta}$. Thus, $\lapp{\app{\rlzSubsetTrans}{t}}{s} \Vdash \fullname{\beta_a} \subseteq \fullname{\delta}$. On the other hand, since $\delta < \beta_a$, $(\fullname{\delta}, \pi) \in \fullname{\beta_a}$ and $\rlzXNotInX \Vdash \fullname{\delta} \not\in \fullname{\delta}$. But this means that $\rlzXNotInX \stackapp \pi \in \falsity{\fullname{\beta_a} \subseteq \fullname{\delta}}$ and hence
    \[
    \lamAbstOne \lamAbstTwo \lapp{\lapp{\app{\rlzSubsetTrans}{\lamtermOne}}{\lamtermTwo}}{\rlzXNotInX} \star t \stackapp s \stackapp \pi \succ \lapp{\app{\rlzSubsetTrans}{t}}{s} \star \rlzXNotInX \stackapp \pi \in \Perp.
    \]
\end{proof}

\begin{lemma}
    $\rlzmodel \Vdash \forall x \exists y (y \rlzin \fullname{\tf{Ord}} \land y \not\in x)$.
\end{lemma}

\begin{proof}
    It will suffice to prove that
    \[
    \lamAbstThree \rapp{\lamtermThree}{\lamAbstOne \lamAbstTwo \lapp{\lapp{\app{\rlzSubsetTrans}{\lamtermOne}}{\lamtermTwo}}{\rlzXNotInX}} \Vdash \forall x (\forall y^{\fullname{\tf{Ord}}} (y \in x) \rightarrow \perp),
    \]
    where $\rlzXNotInX$ and $\rlzSubsetTrans$ are the fixed realizers from the previous lemma.
    
    Fix $a \in \rlzstr$, $\pi \in \Pi$ and $t \Vdash \forall y^{\fullname{\tf{Ord}}} (y \in a)$. By Lemma \ref{lemma: reish ord is proper class}, we can find some ordinal $\beta_a$ such that $\lamAbstOne \lamAbstTwo \lapp{\lapp{\app{\rlzSubsetTrans}{\lamtermOne}}{\lamtermTwo}}{\rlzXNotInX} \Vdash \fullname{\beta_a} \not\in a$. Thus, since $t \Vdash \fullname{\beta_a} \in a$, we have
    \[
    t \star \lamAbstOne \lamAbstTwo \lapp{\lapp{\app{\rlzSubsetTrans}{\lamtermOne}}{\lamtermTwo}}{\rlzXNotInX} \stackapp \pi \in \Perp,
    \]
    as required.
\end{proof}

Next, by using the same argument as in \Cref{theorem:DalethnLargestElement} we can show that $\fullname{\tf{Ord}}$ is closed under initial segments \Rcomment{(in the sense that for all $x \in \fullname{\tf{Ord}}$, if $y \in x$ then $y \in \fullname{\tf{Ord}}$)}, and thus will corresponds to the class of extensional ordinals. 

\begin{lemma}
    $\lamAbstOne \lapp{\app{\lamtermOne}{\rlzXSubsetX}}{\rlzXSubsetX} \Vdash \forall y^{\fullname{\tf{Ord}}} ( y \subseteq \fullname{\tf{Ord}})$, where $\rlzXSubsetX \in \mathcal{R}$ is the fixed realizer such that $\rlzXSubsetX \Vdash \forall x (x \subseteq x)$.
\end{lemma}

\begin{proof}
    Fix $\beta \in \tf{Ord}$ and $t \stackapp \pi \in \falsity{\fullname{\beta} \subseteq \fullname{\tf{Ord}}}$. This means that we can fix some $\delta \in \beta$ such that $t \Vdash \fullname{\delta} \not\in \fullname{\tf{Ord}}$. But, since $\delta$ is an ordinal, 
    \[
    \rlzXSubsetX \stackapp \rlzXSubsetX \stackapp \pi \in \falsity{ \fullname{\delta} \not\in \fullname{\tf{Ord}}} = \bigcup_c \{ s \stackapp s' \stackapp \sigma \divline (c, \sigma) \in \fullname{\tf{Ord}}, \, s s' \Vdash \fullname{\delta} \simeq c \}.
    \]
    Hence $\lamAbstOne \lapp{\app{\lamtermOne}{\rlzXSubsetX}}{\rlzXSubsetX} \star t \stackapp \pi \succ t \star \rlzXSubsetX \stackapp \rlzXSubsetX \stackapp \pi \in \Perp$.
\end{proof}

Thus, $\fullname{\tf{Ord}}$ is realized to be an extensional proper class, every element of which is an $\rlzin$-ordinal, which is closed under initial segments. Since being an $\rlzin$-ordinal implies being an $\in$-ordinal, and in \tf{ZF} the ordinals can be defined as the unique proper class of ordinals which is closed under initial segments, the following is then immediate.

\begin{theorem}
    $\rlzmodel \Vdash \forall x (\tf{Ord}_\in(x) \rightarrow x \in \fullname{\tf{Ord}})$. That is, in any realizability model, $\fullname{\tf{Ord}}$ is a name for the class of extensional ordinals.
\end{theorem}

\subsection{Computing \texorpdfstring{$\omega$}{omega}}

We now proceed to show that $\fullname{\omega}$ corresponds to $\omega$. For this we first consider the lift of the successor function to $\fullname{\tf{Ord}}$ rather than $\hat{\kappa}$. Due to the uniform nature of the realizers used in \Cref{theorem:RecursiveSuccessors}, there is a single realizer $\Rcomment{r} \in \mathcal{R}$ such that for any ordinal $\alpha$,
\[
\Rcomment{r} \Vdash \lquote\fullname{\alpha + 1} \text{ is the successor of } \fullname{\alpha}\rquote.
\]
Hence, the following will be immediate using the general theory of lifting functions (see \cite{Matthews2023}, Section 16 for details). Here we let $\ff{succ}$ denote the successor function on ordinals, $\alpha \mapsto \alpha + 1$.

\begin{lemma}
    Let $\PlusOne = \{ (\op(\fullname{\alpha}, \fullname{\alpha + 1}), \pi) \divline \alpha \in \tf{Ord}, \pi \in \Pi \}$. Then, in $\rlzmodel$, $\PlusOne$ is an $\rlzin$-function from $\fullname{\tf{Ord}}$ to itself. Moreover,
    \[
    \rlzmodel \Vdash \forall x^{\fullname{\tf{Ord}}} (\lquote\PlusOne(x) \text{ is the successor of } x\rquote).
    \]
\end{lemma}

\begin{theorem}
    Let $\rlzXSimeqX$ be a fixed realizer such that $\rlzXSimeqX \Vdash \forall x (x \simeq x)$. Then
    \[
    \lamAbstOne \app{\rlzXSimeqX}{\lamtermOne} \Vdash \forall x^{\fullname{\omega}} (\forall y^{\fullname{\omega}} (x \not\simeq \PlusOne(y)) \rightarrow \forall y (y \notrlzin x)).
    \]
    Hence, $\rlzmodel \Vdash \fullname{\omega} \simeq \boldsymbol{\omega}$, where $\boldsymbol{\omega}$ is extensionally the first infinite ordinal.
\end{theorem}

\begin{proof}
    Fix $n \in \omega$, $t \Vdash \forall y^{\fullname{\omega}} ( \fullname{\Rcomment{n}} \not\simeq \PlusOne(y))$ and $\pi \in \falsity{\forall y (y \notrlzin \fullname{\Rcomment{n}})}$. Since $\falsity{\forall y (y \notrlzin \fullname{n})} = \bigcup_{m < n} \falsity{\fullname{m} \notrlzin \fullname{n}}$, for this to be non-empty we must have that $n \neq 0$. From this we obtain that $t \Vdash \fullname{n} \not\simeq \PlusOne(\fullname{n-1})$. On the other hand, by definition of the lift, this gives $t \Vdash \fullname{n} \not\simeq \fullname{n}$. Thus $\lamAbstOne \app{\rlzXSimeqX}{\lamtermOne} \star t \stackapp \pi \succ \rlzXSimeqX \star t \stackapp \pi \in \Perp$.

    For the second part of the theorem, again using the definition of $\PlusOne$, by taking logical equivalences
    \[
    \rlzmodel \Vdash \forall x \rlzin \fullname{\omega} ( x \simeq \fullname{0} \lor \exists y \rlzin \fullname{\omega} (\lquote x \text{ is the successor of } y\rquote)).
    \]
    Moreover, $\rlzmodel \Vdash 0 \rlzin \fullname{\omega} \land \forall x \rlzin \fullname{\omega} (\PlusOne(x) \rlzin \fullname{\omega})$. From this it follows that $\rlzmodel \Vdash \fullname{\omega} \simeq \boldsymbol{\omega}$ as $\boldsymbol{\omega}$ is the unique set containing $0$ which is closed under successors and satisfies that every element of it is either $0$ or a successor.
\end{proof}

The same method can also be used to show that $\hat{\omega}$ corresponds to $\boldsymbol{\omega}$. As with the successor case, this will not involve any additional special instructions. In particular, by taking $\nu_n = \underline{n}$ in the enumeration of $|\Lambda|$, we can arrange that $\hat{\omega}$ is a well-defined $\rlzin$-TOD which represents $\boldsymbol{\omega}$ in \emph{any} realizability model. We remark here that when $\Lambda$ is countable then the enumeration of terms in this way will have order-type $\omega + \omega$, allowing us to define $\widehat{\omega + \omega}$. However, this makes no difference to the overall analysis.

In order to do this, we first work with the intermediary set $\cjgimel(\ff{dom}(\hat{\omega}))$. Here, the exact same argument will give us that this is a name for $\boldsymbol{\omega}$.

\begin{theorem}
    $\lamAbstOne \app{\rlzXSimeqX}{\lamtermOne} \Vdash \forall x^{\cjgimel(\ff{dom}(\hat{\omega}))} (\forall y^{\cjgimel(\ff{dom}(\hat{\omega}))} (x \not\simeq \widehat{\ff{succ}}(y)) \rightarrow \forall y (y \notrlzin x))$. 
    
    Hence, $\rlzmodel \Vdash \cjgimel(\ff{dom}(\hat{\omega})) \simeq \boldsymbol{\omega}$.
\end{theorem}

\begin{corollary}
    $\rlzmodel \Vdash \hat{\omega} \simeq \boldsymbol{\omega}$.
\end{corollary}

\begin{proof}
    First recall that $\rlzmodel \Vdash \hat{\omega} \subseteq_{\rlzin} \cjgimel(\ff{dom}(\hat{\omega}))$ and therefore, in the extensional model where these are both ordinals, $\hat{\omega} \leq \cjgimel(\ff{dom}(\hat{\omega}))$.

    Next, by \Cref{theorem:HatLiftIsSuccessor}, we have that $\rlzmodel \Vdash \forall x^{\hat{\kappa}} (\lquote \widehat{\ff{succ}}(x) \text{ is the successor of } x\rquote$). Taking the natural restriction of this map gives us an $\rlzin$-function from $\hat{\omega}$ to $\cjgimel(\ff{dom}(\hat{\omega}))$. Moreover, by the uniform nature of the realizers from \Cref{theorem:HatLiftProperties} and \Cref{theorem:HatLiftIsSuccessor}, we have
    \[
    \rlzmodel \Vdash \forall x^{\hat{\omega}} ( x \in \widehat{\ff{succ}}(x) \land \widehat{\ff{succ}}(x) \rlzin \hat{\omega}).
    \]
    Thus, $\hat{\omega}$ is realized to be an infinite ordinal which is at most $\boldsymbol{\omega}$ and therefore the two must be extensionally equal.
\end{proof}

The next thing to show would be that \Lcomment{$\fullname{\alpha}$ always corresponds to the $\alpha$-th ordinal of the realizability model}. \Rcomment{This is indeed true in the case for forcing where one shows that the interpretation of the canonical name for $\alpha$, denoted $\check{\alpha}^G$, is always equal to $\alpha$.} Note that, for $\delta =  (\omega_1)^{\tf{V}}$, it need not be the case that $\fullname{\delta}$ is the $\omega_1$ of $\rlzmodel$, for example it may be collapsed \Lcomment{to a countable ordinal, and that could happen even in forcing models}. Therefore we must detach an ordinal from any cardinal it could represent in a given model. 
\Rcomment{However, proving a general correspondence between $\fullname{\alpha}$ and $\alpha$ appears difficult to achieve. For example, it was relatively straightforward to show that $\fullname{\omega}$ was an infinite ordinal, and thus $\rlzmodel \Vdash \boldsymbol{\omega} \leq \fullname{\omega}$. The difficulty is in the reverse direction. Here we heavily used the fact that $\boldsymbol{\omega}$ is the unique non-empty set which is closed under successors and satisfies that every element is either $0$ or a successor. Without such a clear description of an ordinal the authors do not know how to prove $\fullname{\alpha}$} \Lcomment{remains the $\alpha$-th ordinal in the realizability model.} 

\section{Preserving Cardinals}

Since $\fullname{\omega}$ is a name for $\omega$, it will remain an $\in$-cardinal in any realizability model. The next question to address is which other ordinals give rise to cardinals. 
As with forcing, in general it \Lcomment{may} not be the case that $\delta$ being a cardinal in \tf{V} implies that $\fullname{\delta}$ should remain a cardinal in $\rlzmodel$, in fact, it is very difficult to give general conditions for when this occurs.

However, in this section we shall show that if $\delta > | \Lambda|$ is a regular cardinal then $\fullname{\delta}$ remains a cardinal in $\rlzmodel$. This is the realizability analogue of the fact that forcing preserves cardinals greater than the size of the poset.

Recall that $\ff{ExtFun}(f, a)$ is the formula expressing that $f$ is an extensional function with domain $a$, namely
\[
\ff{ExtFun}(f, a) \equiv \forall x_1^a x_2^a, y_1, y_2 \, (x_1 \simeq x_2 \rightarrow \op(x_1, y_1) \rlzin f \rightarrow \op(x_2, y_2) \rlzin f \rightarrow y_1 \not\simeq y_2 \rightarrow \perp).
\]

\begin{theorem} \label{theorem:preservingcardinals}
Let $\delta > |\Lambda|$ be a regular cardinal. Then 
\[
\rlzmodel \Vdash \forall f \, \forall a \rlzin \fullname{\delta} \, \exists b \rlzin \fullname{\delta} (\ff{ExtFun}(f, a) \rightarrow \forall y \rlzin a (\op(y, b) \rlzin f \rightarrow \perp)).
\]
Namely, in $\rlzmodel$, if $a \rlzin \fullname{\delta}$ then any extensional function $f \colon a \rightarrow \fullname{\delta}$ is not an $\rlzin$-surjection. Thus, $\fullname{\delta}$ is an $\rlzin$-cardinal in $\rlzmodel$ and hence also an $\in$-cardinal.
\end{theorem}

\begin{proof}
Fix $\rlzXSimeqX, \rlzInIsNotSimeq \in \mathcal{R}$ to be realizers such that $\rlzXSimeqX \Vdash \forall x (x \simeq x)$ and $\rlzInIsNotSimeq \Vdash \forall x \forall y (x \rlzin y \rightarrow x \not\simeq y)$ and let 
\[
u_0 \coloneqq \lamAbstTwo \lamAbstThree \lamtermTwo \rlzXSimeqX \lamtermThree \lamtermThree (\app{\rlzInIsNotSimeq}{\identity}).
\]
We shall first show that for any $f \in \rlzstr$ and $\alpha \in \delta$ there is some $\mu \in \delta$ such that $u_0 \Vdash \varphi(\fullname{\alpha}, \fullname{\mu}, f)$ where 
\[
\varphi(a, b, f) \equiv \ff{ExtFun}(f, a) \rightarrow \forall y^a (\op(y, b) \rlzin f \rightarrow \perp).
\]
To do this, for $\gamma \in \alpha$ and $s \in \Lambda$, let $X_{\gamma, s} \coloneqq \{ \beta \in \delta \divline s \Vdash \op(\fullname{\gamma}, \fullname{\beta}) \rlzin f \}$. Since $|\Lambda|$ and $\alpha$ are both less than $\delta$, which is a regular cardinal, 
\[
\mu \coloneqq \ff{sup}( \{ \ff{min}(X_{\gamma, s}) \divline s \in \Lambda, \gamma \in \alpha \}) + 1 < \delta.
\]
Note that if $\beta < \mu$ then $\identity \Vdash \fullname{\beta} \rlzin \fullname{\mu}$ and therefore $\app{\rlzInIsNotSimeq}{\identity} \Vdash \fullname{\beta} \not\simeq \fullname{\mu}$.

Now, suppose that $t \Vdash \ff{ExtFun}(f, a)$ and $r \stackapp \pi \in \falsity{\forall y^{\fullname{\alpha}} (\op(y, \fullname{\mu}) \rlzin f \rightarrow \perp)}$. Then we can fix $\gamma \in \alpha$ such that $r \Vdash \op(\fullname{\gamma}, \fullname{\mu}) \rlzin f$. Thus we must have that $X_{\gamma, r} \neq \emptyset$, so let $\beta = \ff{min}(X_{\gamma, r})$. Then, $r \Vdash \op(\fullname{\gamma}, \fullname{\beta}) \rlzin f$ and, since $\beta < \mu$, $\app{\rlzInIsNotSimeq}{\identity} \Vdash \fullname{\beta} \not\simeq \fullname{\mu}$. From this we can conclude that $u_0 \star t \stackapp r \stackapp \pi \succ t \stackapp \rlzXSimeqX \stackapp r \stackapp r \stackapp \app{\rlzInIsNotSimeq}{\identity} \stackapp \pi \in \Perp$. \\

\noindent To finish the proof, it will suffice to show that
\[
u \coloneqq \lamAbst{k} \app{k}{u_0} \Vdash \forall f \, \forall a^{\fullname{\delta}} \, (\forall b^{\fullname{\delta}} (\varphi(a, b, f) \rightarrow \perp) \rightarrow \perp).
\]
Fix $f \in \rlzstr$, $\alpha \in \delta$, $\pi \in \Pi$ and let us suppose that $t \Vdash \forall b^{\fullname{\delta}} (\varphi(\fullname{\alpha}, b, f) \rightarrow \perp)$. Taking $\mu$ as in the claim, we have $u_0 \Vdash \varphi(\fullname{\alpha}, \fullname{\mu}, f)$. On the other hand, by assumption, $t \Vdash \varphi(\fullname{\alpha}, \fullname{\mu}, f) \rightarrow \perp$, and thus 
\[
\lamAbst{k} \app{k}{u_0} \star t \stackapp \pi \succ t \star u_0 \stackapp \pi \in \Perp,
\]
completing the proof
\end{proof}

The same argument in fact shows us that in $\rlzmodel$, for every $a$ there is an $\rlzin$-cardinal $\delta$ such that $a$ does not surject onto $\delta$.

\begin{theorem}
Given $a \in \rlzstr$, let $\delta = | \ff{dom}(a) \cup \Lambda|^+$. Then, in $\rlzmodel$, $\fullname{\delta}$ is an $\rlzin$-cardinal and there is no extensional function which is a $\rlzin$-surjection of $a$ onto $\fullname{\delta}$.
\end{theorem}

\begin{proof}
The proof is very similar to \Cref{theorem:preservingcardinals} and therefore we will just sketch the details. Since $\identity \Vdash \forall z (z \notrlzin \cjgimel(\ff{dom}(a)) \rightarrow z \notrlzin a)$, we have $\rlzmodel\Vdash a \subseteq \cjgimel(\ff{dom}(a))$. So, it will suffice to prove there is no surjection of $\cjgimel(\ff{dom}(a))$ onto $\fullname{\delta}$. Fix $a \in \rlzstr$ and let $\delta = | \ff{dom}(a) \cup \Lambda|^+$. Then, by \Cref{theorem:preservingcardinals}, $\fullname{\delta}$ is a cardinal in $\rlzmodel$. Fix $f \in \rlzstr$ and, for $b \in \ff{dom}(a)$ and $s \in \Lambda$, let $X_{b, s} \coloneqq \{ \beta \divline s \Vdash \op(b, \fullname{\beta}) \rlzin f \}$. Since $|\ff{dom}(a)|$ and $|\Lambda|$ are both less than $\delta$, we can fix
\[
\mu = \ff{sup}(\{ \ff{min}(X_{b, s}) \divline s \in \Lambda, \beta \in \alpha \}) + 1 < \delta.
\]
Then, proceeding as before, if $t \Vdash \ff{ExtFun}(f, \cjgimel(\ff{dom}(a)))$ then $\app{u_0}{t} \Vdash \forall b^a \op(b, \fullname{\mu}) \notrlzin f$. The end of the argument from \Cref{theorem:preservingcardinals} can then be used to show that the same realizer $u$ will realize that there is no extensional function which is a $\rlzin$-surjection of $\cjgimel(\ff{dom}(a))$ onto $\fullname{\delta}$.
\end{proof}

These arguments do not immediately give us that $\fullname{\delta}$ remains a \emph{regular} cardinal. However, a very slight modification of the proof can be used to show that for every $a \rlzin \fullname{\delta}$, every $\rlzin$-function $f \colon a \rightarrow \fullname{\delta}$ is \emph{bounded}. Thus, \Rcomment{by \Cref{theorem:NEACGivesCardinals},} in the presence of \tf{NEAC} $\fullname{\delta}$ is indeed a regular $\in$-cardinal.

\begin{corollary}\label{theorem:FullnameKappaPlusCardinal}
    Let $\delta > |\Lambda|$ be a regular cardinal. Then
    \[
    \rlzmodel \Vdash \forall f \forall a^{\fullname{\delta}} \exists b^{\fullname{\delta}} (\rlzin\ff{-Fun}(f, a) \rightarrow \forall y^{a} \forall z^{\fullname{\delta}} (\ff{op}(y, z) \rlzin f \rightarrow z \rlzin b)).
    \]
    That is, for every $a \rlzin \fullname{\delta}$, every $\rlzin$-function $f \colon a \rightarrow \fullname{\delta}$ is bounded.

    Hence, if $\rlzmodel \Vdash \tf{NEAC}$, then $\fullname{\delta}$ is a regular $\in$-cardinal in $\rlzmodel$.
\end{corollary}

Using the argument of \Cref{theorem:preservingcardinals} we can in fact show that \emph{every} cardinal greater than the size of the realizability algebra is preserved.

\begin{theorem}
    Let $\delta > |\Lambda|$ be an infinite cardinal. Then, in $\rlzmodel$, $\fullname{\delta}$ is an $\in$-cardinal. 
\end{theorem}

\begin{proof}
    It only remains to consider the case when $\delta$ is a singular cardinal. Then, since Choice holds in the ground model, $\delta$ is the supremum of a set, $S$, of regular cardinals, all of which are greater than $|\Lambda|$. Let $\lift{S} = \{ (\fullname{\alpha}, \pi) \divline \alpha \in S, \pi \in \Pi \}$. Since $\lift{S}$ can be written as a $\cjgimel$-name, \Cref{theorem:BoundedUniversals} gives us that $\forall y^{\lift{S}}$ corresponds to the bounded universal quantifier.
    Then, by the uniform nature of proof of \Cref{theorem:preservingcardinals}, in particular that it did not depend on $\delta$,
    \[
    \rlzmodel \Vdash \forall y^{\lift{S}} \; \lquote y \text{ is an} \in\text{-cardinal}\rquote.
    \]
    Next, it is easy to see that
    \[
    \identity \Vdash \forall y^{\lift{S}} y \rlzin \fullname{\delta}
    \quad \text{and} \quad
    \identity \Vdash \forall x^{\fullname{\delta}} (\forall y^{\lift{S}} (x \notrlzin y) \rightarrow \perp).
    \]
    Using the fact that $(\rlzstr, \in, \simeq)$ gives a model of \tf{ZF}, it follows that
    \[
    \rlzmodel \Vdash \bigcup \lift{S} \simeq \fullname{\delta},
    \]
    where $\bigcup$ is the internal union in the extensional part of the realizability model. Thus, $\fullname{\delta}$ is the supremum of a set of $\in$-cardinals, and hence also an $\in$-cardinal.
\end{proof}

\section{Epsilon Totally Ordered Sets}

While the class $\fullname{\tf{Ord}}$ extensionally corresponds to the class of ordinals in a realizability model, we again remark that the definition of an ordinal as an $\rlzin$-transitive set of $\rlzin$-transitive sets is not always the most appropriate. Therefore, it is worth considering what is the class of $\rlzin$-TODs. Since, under the reasonable assumption of having the special instruction $\chi$, $\hat{\alpha}$ is an $\rlzin$-TOD for all $\alpha < \kappa \coloneqq | \Lambda |$, $\hat{\kappa}$ will in general provide a lower bound for the supremum of the class of $\rlzin$-TODs. We next show that every such ordinal in fact injects into $\fullname{\kappa^+}$ and thus forms a set. In order to do that, we need to first investigate certain classes which can be built using elements of $\fullname{2}$.

\subsection{Building Classes from \texorpdfstring{$\fullname{2}$}{recursive 2}}

\Rcomment{Recall that $2$ can be seen as a Boolean algebra with binary operations $\lor$ and $\land$ and unary operation $\neg$. Lifting these operations gives a new Boolean algebra on the set of truth values $\fullname{2}$, the proof of which can be found in Section 16 of \cite{Matthews2023}.

\begin{theorem}
    $\fullname{2}$ is a Boolean algebra with binary operations $\lor$ and $\land$ and unary operation $\neg$ given by the lift of the corresponding (trivial) operations on $\{\fullname{0}, \fullname{1}\}$.
\end{theorem}

We remark here that this Boolean algebra may satisfy many interesting properties. For example, it may be \emph{atomless} ($a$ is said to be an \emph{atom} if, for any $x$, $a \land x$ is either $0$ or $a$) a property which is used in \cite{Krivine2012} to produce a realizability model in which the reals are not well-ordered.
}

We \Rcomment{now} define a function $\mathfrak{h} \colon \Rcomment{\{ \fullname{0}, \fullname{1} \}} \times \rlzstr \rightarrow \rlzstr$, by
\[
\mathfrak{h}(\Rcomment{\fullname{a}},x) \mapsto 
\begin{cases}
    0 & \text{ if } a = 0, \\
    x & \text{ if } a = 1.
\end{cases}
\]
\Rcomment{Since $\{ \fullname{0}, \fullname{1} \} = \ff{dom}(\fullname{2})$, by \Cref{definition:FunctionLift} this} lifts to \Rcomment{an $\rlzin$-}function $\lift{\mathfrak{h}} \colon \fullname{2} \times \rlzstr \rightarrow \rlzstr$ which we will also denote by $ax$.

Two important properties of $\lift{\mathfrak{h}}$ are \Rcomment{``\emph{linearity}'', as defined in \cite{Krivine2015}, and that} if two elements $a, b$ of $\fullname{2}$ are incompatible then \Rcomment{$abx = \lift{\mathfrak{h}}(a, \lift{\mathfrak{h}}(b, x))$ results in} a set which is equal to $0.$

\begin{lemma}[Linearity] \label{theorem:LinearityOnFullname2}
    Let \Rcomment{$f \in \rlzstr$ be a binary function satisfying any of the concepts of being a function from \Cref{Section:Functions}}. Then
    \[
    \identity \Vdash \forall a^{\fullname{2}} \forall b^{\fullname{2}} \forall x \forall y (a b f(x, y) = a b f(ax, by)).
    \]
\end{lemma}

\begin{proof}
    It suffices to prove that for any $x, y \in \rlzstr$, $\identity \Vdash a b f(x, y) = a b f(ax, ay)$ for the 4 possible values of the pair $(a, b) \in \{0, 1\}^2$. But this is just $\identity \Vdash \fullname{0} = \fullname{0}$, unless $a  = b = 1$ in which case we have $\identity \Vdash f(x, y) = f(x, y)$.
\end{proof}

\begin{proposition} \label{theorem:anotaxIs0}
    $\rlzmodel \Vdash \forall x \forall a^{\fullname{2}} \forall b^{\fullname{2}} (a \land b = \fullname{0} \rightarrow abx = \fullname{0})$.
\end{proposition}

On the other hand, it is worth remarking that, from the point of view of the realizability structure, if $a \neq \Rcomment{\fullname{1}}$ then $\lift{\mathfrak{h}}(a, x)$ is always extensionally equal to the empty set.

\begin{proposition}
    $\identity \Vdash \forall x \forall a^{\fullname{2}} ( a \neq \fullname{1} \rightarrow \forall y (y \notrlzin ax))$.
\end{proposition}

\begin{proof}
    Fix $x \in \rlzstr$, there are two cases to consider. If $a = 0$ then $\falsity{\fullname{0} \neq \fullname{1} \rightarrow \forall y (y \notrlzin \fullname{0})} = \emptyset$ and therefore any term realizes the statement. On the other hand, if $t \Vdash \fullname{1} \neq \fullname{1}$ then for any $\pi \in \Pi$, $t \star \pi \in \Perp$.
\end{proof}

We next define a ground model function $\langle \cdot \Rcomment{\; < \;} \cdot \rangle \colon \rlzstr^2 \rightarrow 2$ by
\[
\langle x < y \rangle = 
\begin{cases}
    1 & \text{ if } x \in \ff{dom}(y) \\
    0 & \text{otherwise}
\end{cases}
\]
This then lifts to a function $\langle \cdot \Rcomment{\; < \;} \cdot \rangle \colon \rlzstr^2 \rightarrow \fullname{2}$ in $\rlzmodel$ which translates the formula $x \rlzin y$ into a function into $2$, that is a truth value.

\begin{proposition} \label{theorem:RlzinPreserves<}
    $\rlzmodel \Vdash \forall x \forall y ( x \rlzin y \rightarrow \langle x < y \rangle = \fullname{1})$.
\end{proposition}

\begin{proof}
    It suffices to prove that for any $a, b \in \rlzstr$, $\identity \Vdash \langle a < b \rangle \neq \fullname{1} \rightarrow a \notrlzin b$, which is immediate from the definition.
\end{proof}

Finally, in $\rlzmodel$, given $a \rlzin \fullname{2}$ we shall define classes $M_a$ such that for any $x$, $x \rlzin M_a$ if and only if there is some $y \rlzin \rlzstr$ such that $x = ay$. Naively, we should have $M_a = \{ x \divline \exists y (x = ay)\}$. The issue is that the symbol $=$ is only defined for sets and the Leibniz definition of equality ($x = y \leftrightarrow \forall z (x \rlzin z \leftrightarrow y \rlzin z)$) does not make sense for proper classes; thus the above specification does not lead to a unique definition for the class $M_a$.

Instead, to ensure that this is a well-defined definition, we formally work in the second-order realizability theory $\GBepsilon$ which is defined in Section 9 of \cite{FontanellaGeoffroyMatthews2024}. Here we will have a two-sorted first-order logic where the first sort will represent sets (denoted by lower case letters) and the second sort classes (denoted by upper case letters). It is then shown that if $(\tf{V}, \mathcal{C})$ is a model of \tf{GB} then we can construct the realizability model $\rlzmodel = (\rlzstr, \mathcal{D})$ as a natural generalisation of the set case and prove that it is a model of $\GBepsilon$.

However, since the details of the theory and extended realizability interpretation are not vital, we will omit all details of the construction itself.

\begin{proposition}\label{prop: Ma}
    $\rlzmodel \Vdash \forall a^{\fullname{2}} \exists M \forall x (x \rlzin M \longleftrightarrow \exists y ( \lift{\mathfrak{h}}(a, y) = x)$.
\end{proposition}

\begin{proof}
    Given $a \in \{0, 1\}$, let $Z_a = \{ (x, t \stackapp \pi) \divline x \in \rlzstr, \, \pi \in \Pi, \, t \Vdash \exists y ( \lift{\mathfrak{h}}(a, y) = x ) \}$. Then, using the proof of Separation in realizability models, for any $x \in \rlzstr$,
    \[
    \falsity{x \notrlzin Z_a} = \falsity{\exists y (\lift{\mathfrak{h}}(a, y) = x) \rightarrow \perp}
    \]
    and thus
    \[
    \identity \Vdash x \notrlzin Z_a \rightarrow (\exists y (\lift{\mathfrak{h}}(a, y) = x) \rightarrow \perp).
    \]
    and
    \[
    \identity \Vdash (\exists y (\lift{\mathfrak{h}}(a, y) = x) \rightarrow \perp) \rightarrow x \notrlzin Z_a.
    \]
    Hence, 
    \[
    \lamAbstOne \lapp{\app{\lamtermOne}{\identity}}{\identity} \Vdash \exists M \forall x (x \notrlzin M \longleftrightarrow (\exists y (\lift{\mathfrak{h}}(a, y) = x) \rightarrow \perp)).
    \]
    Since the same realizer works for both $a = 0$ and $a = 1$, it follows that $\lamAbstOne \lapp{\app{\lamtermOne}{\identity}}{\identity}$ can be used to construct a realizer for our desired statement.
\end{proof}

\begin{remark} 
We call $M_a$ any class satisfying the statement of Proposition \ref{prop: Ma}. 
Observe that there may be more than one class satisfying this definition since $M_a$ is only defined up to $\rlzin$-\Rcomment{equality} (having the same $\rlzin$-elements), which is strictly weaker than the Leibniz equality given by $=$.

\end{remark}

We shall denote this weaker sense of equality as follows, where we use the $\simeq$ symbol to remind the reader that it is closer in flavour to extensional equality.

\begin{definition}
    Say that $x \simeq_{\rlzin} y$ iff $\forall z (z \rlzin x \longleftrightarrow z \rlzin y)$.
\end{definition}

\subsection{The Set of all Epsilon Totally Ordered Sets}

We next show that the collection of all $\rlzin$-ToDs form a set. To do this, we first need a lemma which provides a method to approximate functions in the realizability model by functions in the ground model.

\begin{lemma}[Approximating Functions using the Ground Model]\label{theorem:GroundModelApprox}
    Let $\mathcal{A}$ be a realizability model and let $\kappa = | \Lambda|$. Then for any formula $\varphi(u, x_1, . . . , x_n)$ in $Fml_{\rlzin}$ there exists a class function $g_\varphi : \kappa \times \rlzstr^n \rightarrow \rlzstr$ such that
    \[
    \rlzmodel \Vdash \forall x_1, \dots, x_n \, \big( \exists u \, \varphi(u, x_1, \dots, x_n) \rightarrow \exists a^{\fullname{\kappa}} \varphi(\lift{g}_\varphi(a, x_1, \dots, x_n), x_1, \dots, x_n) \big),
    \]
    where $\lift{g}_\varphi$ is the canonical lift of $g_\varphi$ to an $\rlzin$-class function.
\end{lemma}

\begin{proof}
    Fix an enumeration $\langle \nu_\alpha \divline \alpha \in \kappa \rangle$ of $\Lambda$ and let $\varphi(u, x)$ be a formula in $Fml_{\rlzin}$. To ease notation, without loss of generality, we shall assume that $\varphi$ has only two free variables. Given $\alpha \in \kappa$, let $P_\alpha = \{ \pi \in \Pi \divline \nu_\alpha \star \pi \not\in \Perp \}$. Using \tf{AC} in the ground model, define $g_\varphi \colon \kappa \times \rlzstr \rightarrow \rlzstr$ such that for any $b \in \rlzstr$ and $\alpha \in \kappa$, if $P_\alpha \cap \falsity{\forall u \, \neg \varphi(u, b)} \neq \emptyset$ then $P_\alpha \cap \falsity{\neg \varphi(g_\varphi(\alpha, b), b} \neq \emptyset$. Otherwise, set $g_\varphi(\alpha, b) = 0$.

    Then the lift of $g_\varphi$ is $\lift{g}_\varphi = \{ ( \op ( \op(\fullname{\alpha}, b), g_\varphi(\alpha, b) ), \pi) \divline b \in \rlzstr, \alpha \in \kappa, \pi \in \Pi \}$. We shall conclude by showing that
    \[
    \identity \Vdash \forall x (\forall a^{\fullname{\kappa}} \neg \varphi(\lift{g}_\varphi(a, x), x) \rightarrow \forall u \neg \varphi(u, x)).
    \]
    For suppose this was not the case. Fix $b \in \rlzstr$, $\nu_\gamma \Vdash \forall a^{\fullname{\kappa}} \neg \varphi(\lift{g}_\varphi(a, b), b)$ and $\pi \in \falsity{ \forall u \neg \varphi(u, b)}$ such that $\identity \star \nu_\gamma \stackapp \pi \not\in \Perp$. From this it follows that $\nu_\gamma \star \pi \not\in \Perp$. Therefore, $\pi \in P_\gamma \cap \falsity{\forall u \neg \varphi(u, b)}$. Now, by the construction of $g_\varphi$, fix $\sigma \in P_\gamma \cap \falsity{\neg \varphi(g_\varphi(\gamma, b), b)}$. But then $\nu_\gamma \Vdash \neg \varphi( \lift{g}_\varphi(\fullname{\gamma}, b), b)$. Moreover, since $\sigma \in P_\gamma$, $\nu_\gamma \star \sigma \not\in \Perp$. 

    Finally, since $\sigma \in \falsity{\neg \varphi(g_\varphi(\gamma, b), b)}$, $\sigma \in \falsity{\forall a^{\fullname{\kappa}} \neg \varphi(\lift{g}_\varphi(a, b), b)}$ and so $\nu_\gamma \star \sigma \in \Perp$, yielding our desired contradiction.
\end{proof}

\begin{lemma}\label{theorem:MaMbNotBothUnbounded} 
Let $\kappa = | \Lambda|$ and suppose that $\rlzmodel \Vdash \tf{NEAC}$ and there exist $a, b \rlzin \fullname{2}$ such that $a, b \neq 0$ and $a \wedge b = 0$. Then, if there are two $\rlzin$-functions $f_0: M_a\to \fullname{Ord}$ and $f_1: M_{b}\to \fullname{Ord},$ then one of them is bounded in $\fullname{\kappa^+}$.
\end{lemma}

\begin{proof}
    Suppose for a contradiction that both $f_0$ and $f_1$ are unbounded in $\fullname{\kappa^+}$. This means that
    \[
    \rlzmodel \Vdash \forall z \rlzin \fullname{\kappa^+} \, \exists x \rlzin M_a \, \exists y \rlzin M_{b} \, z \rlzin f_0(x) \cap f_1(y).
    \]
    Extend $f_0$ and $f_1$ to $\rlzstr$ by setting $f_0(x) = \fullname{0}$ if $x \notrlzin M_a$ and $f_1(x) = \fullname{0}$ if $x \notrlzin M_{b}$. Let $\psi(x, y) \equiv f_0(x) \rlzin f_1(by)$.

    Then, by \Cref{theorem:GroundModelApprox}, there is a ground model function $g$ which lifts to an $\rlzin$-function $\lift{g}$ such that 
    \[
    \forall x \big( \exists y \, \psi(x, y) \rightarrow \exists \sigma^{\fullname{\kappa}} \psi(x, \lift{g}(x, \sigma)) \big).
    \]
    In particular, 
    \[
    \big( \exists y (f_0(ax) \rlzin f_1(by)) \rightarrow \exists \sigma^{\fullname{\kappa}} (f_0(ax) \rlzin f_1(b \lift{g}(ax, \sigma)) ) \big).
    \]
    Now, since $a \wedge b = 0$, $b \lift{g}(ax, \sigma) = b \lift{g}(bax, \sigma) = b \lift{g}(0, \sigma)$. So
    \[
    \big( \exists y (f_0(ax) \rlzin f_1(by)) \rightarrow \exists \sigma^{\fullname{\kappa}} (f_0(ax) \rlzin f_1(b \lift{g}(0, \sigma)) ) \big).
    \]

    Next, by hypothesis, $\forall z \rlzin \fullname{\kappa^+} \, \exists x, y (z \rlzin f_0(ax) \rlzin f_1(by))$. Thus,
    \[
    \forall z \rlzin \fullname{\kappa^+} \, \exists x \, \exists \sigma^{\fullname{\kappa}} ( z \rlzin f_0(ax) \rlzin f_1(b\lift{g}(ax, \sigma))). 
    \]
     So, using the fact that any element of $\fullname{\kappa^+}$ is an $\rlzin$-transitive set and $b\lift{g}(ax, \sigma) = b\lift{g}(0, \sigma)$, $\forall z \rlzin \fullname{\kappa^+} \, \exists \sigma^{\fullname{\kappa}} (z \rlzin f_1(b\lift{g}(0, \sigma))$. Using \tf{NEAC}, we can define a function $h \colon \fullname{\kappa^+} \rightarrow \fullname{\kappa}$ such that $h(z) = \sigma_z$ for some $\sigma_z$ such that $z \rlzin f_1(b\lift{g}(0, \sigma_z))$. But, from this it follows that the function $h' \colon \fullname{\kappa} \rightarrow \fullname{\kappa^+}$, $\sigma \mapsto f_1(b\lift{g}(0, \sigma))$ is unbounded in $\fullname{\kappa^+}$, contradicting \Cref{theorem:FullnameKappaPlusCardinal}.
\end{proof}

\begin{theorem} Suppose $\fullname{2}$ is non-trivial and $\rlzmodel \Vdash \tf{NEAC}$, then any $\rlzin$-function from an $\rlzin$-TOD into $\fullname{\kappa^+}$ is bounded.
\end{theorem}

\begin{proof}
    Let $X$ be an $\rlzin$-TOD. First, observe that for any $a \rlzin \fullname{2}$ with $a \neq \fullname{0}$, the function $x \mapsto ax$ defines an $\rlzin$-injection from $X$ into $M_a$. Thus is because, if $x \neq y$ in $X$ then either $x \rlzin y$ or $y \rlzin x$. Without loss of generality, suppose that $x \rlzin y$. Then, by \Cref{theorem:RlzinPreserves<}, $\langle x < y \rangle = \fullname{1}$ and therefore $\langle ax < ay \rangle = a$. Thus $ax \neq ay$ since $a \neq \fullname{0}$.

    Now, suppose that $f \colon X \rightarrow \fullname{\kappa^+}$ were an unbounded $\rlzin$-function. By the above observation, there is an $\rlzin$-surjection $g \colon M_a \rightarrow X$ for all $a \rlzin \fullname{2}$, $a \neq \fullname{0}$. But then $g \circ f \colon M_a \rightarrow \fullname{\kappa^+}$ must be an unbounded $\rlzin$-function.

    But now, if $\fullname{2}$ is non-trivial then we can fix some $a \rlzin \fullname{2}$ such that $a, \lnot a \neq \fullname{0}$. Then we must have unbounded $\rlzin$-functions from $M_a$ and $M_{\lnot a}$ onto $\fullname{\kappa^+}$, contradicting \Cref{theorem:MaMbNotBothUnbounded}. 
\end{proof}

\begin{theorem}
    Suppose $\fullname{2}$ is non-trivial and \tf{NEAC} holds in $\rlzmodel$. Then $\{ X \divline X \text{ is a } \rlzin\text{-TOD} \} \subseteq \fullname{\kappa^+}$. Hence, the collection of all $\rlzin$-TODs extensionally forms a set.
\end{theorem}

\begin{proof}
    First note that any $\rlzin$-TOD $X$ is automatically an $\rlzin$-ordinal and hence an $\in$-ordinal. Next, since any $\rlzin$-function from $X$ into $\fullname{\kappa^+}$ is bounded, it must also be the case that every $\in$-function from $X$ into $\fullname{\kappa^+}$ is bounded. Thus, since $X$ and $\fullname{\kappa^+}$ are both $\in$-ordinals, they can be compared, and it must be $X \in \fullname{\kappa^+}$. 
\end{proof}

\begin{remark}
    If $\fullname{2}$ is non-trivial then in general the collection of $\rlzin$-TODs does not form a set in the non-extensional sense. This is because, for any $a \rlzin \fullname{2}$ with $a \neq \fullname{1}$, $ax \simeq_{\rlzin} \fullname{0}$ for every $x$. Thus there are no $\rlzin$-elements of $ax$ and so $ax$ trivially satisfies the definition of being an $\rlzin$-TOD. 

    On the other hand, if $\fullname{2}$ is trivial then it can be shown that $\fullname{\alpha}$ is an $\rlzin$-TOD for every ordinal $\alpha$, and thus they form a proper class. 
\end{remark}

\section{Chain Conditions}\label{Section:ChainConditions}

In this section we study a version of the $\delta$-chain condition which has been adapted to work with realizability algebras. Recall that the $\delta$-chain condition for forcing asserts that every antichain, that is a set of incompatible conditions, has cardinality strictly less than $\delta$. For Boolean algebras, being incompatible is a symmetric property ($p$ is incompatible with $q$ if and only if $q$ is incompatible with $p$). However, in our system, we will need an ordered version of the notion of antichain. This is because of the non-symmetric way we realize conjunctions, namely $t \Vdash a \simeq b$ does not in general imply that $t \Vdash b \simeq a$. \Lcomment{Of course the axioms of equality are realized, hence $\mathcal{N}\models a\simeq b \imp b\simeq a,$ but one should just be aware of the fact that realizing $a\simeq b$ and $b\simeq a$ may require two distinct realizers.} 

\begin{definition}
Let us say that a realizability algebra satisfies the \emph{$\delta$-chain condition} if there exists a realizer $\rlzfont{p} \in \rlzset$ such that for every $\delta$-sequence of terms $\langle u_\beta \divline \beta \in \delta \rangle \subseteq \Lambda$, for every $t \in \Lambda$ and $\pi \in \Pi$:
\[
\textrm{if } \forall \gamma, \beta \in \delta (\gamma < \beta \rightarrow (t \star u_\gamma \stackapp u_\beta \stackapp \pi \in \Perp)), \; \textrm{ then } \exists \beta \in \delta \, \rlzfont{p} \star t \stackapp u_\beta \stackapp \pi \in \Perp.
\]
\end{definition}

We also have a \Rcomment{weaker, but more natural,} uniform version of the Chain Condition which is more symmetric in nature.

\begin{definition}
Let us say that a realizability algebra satisfies the \emph{uniform $\delta$-chain condition} if there exists a realizer $\rlzfont{p} \in \rlzset$ such that for every $A \subseteq \Lambda$ of cardinality at least $\delta$, for every $t \in \Lambda$ and $\pi \in \Pi$:
\[
\text{if for every } a \neq b \text{ in } A \; (t \star a \stackapp b \stackapp \pi \in \Perp), \; \text{ then there exists an } a \in A \text{ such that } \rlzfont{p} \star t \stackapp a \stackapp \pi \in \Perp.
\]
\end{definition}

To provide evidence that this is an appropriate way to define the realizability analogue of chain conditions, we begin by proving that if $\mathbb{B}$ is a Boolean algebra which satisfies the $\delta$-chain condition then it induces a realizability algebra that also satisfies the $\delta$-chain condition. This algebra is formally defined in Section 2.2 of \cite{FontanellaGeoffroy2020} and its properties are explicitly studied in Section 19 of \cite{Matthews2023}. For completeness, we sketch the construction now.

Given a Boolean algebra $\mathbb{B} = (B, 1, 0, \land, \lor, \neg)$ \Rcomment{we define the realizability algebra $\mathcal{A}_{\mathbb{B}}$ as follows:} $\mathcal{A}_{\mathbb{B}} = (\Rcomment{\emptyset, \{ \omega_b \divline b \in B \}}, \prec, \Rcomment{\Perp})$, so \Rcomment{there are no special instructions and} there is a stack bottom for every condition of $\mathbb{B}$. To define the pre-order and the pole, we begin by inductively defining a function $\tau \colon \Lambda^{\star}_{(0, \mu)} \cup \Pi_{(0, \mu)} \rightarrow \mathbb{B}$, where $\Lambda^{\star}_{(0, \mu)}$ denotes the set of all - possibly open - $\lambda_c$-terms.
\begin{itemize}
    \item for every stack bottom $p$, we let $\tau(p) \coloneqq p$;
    \item for every term $t$ and every stack $\pi$, we let $\tau(t \star \pi) \coloneqq \tau(t) \land \tau(\pi)$;
    \item for every variable $x$, $\tau(x) \coloneqq \tau(\cc) \coloneqq 1$;
    \item for all $\lambda_c$-terms $t, u$, we let $\tau(tu) \coloneqq \tau(t) \land \tau(u)$;
    \item for every variable $x$ and every term $t$, we let $\tau(\lambda x \lambdaapp t) \coloneqq \tau(t)$;
    \item for every stack $\pi$, we let $\tau(\saverlz{\pi}) \coloneqq \tau(\pi)$.
\end{itemize}
Then, we define $\prec$ by $t \star \pi \succ s \star \sigma$ if and only if $\tau(t) \land \tau(\pi) \leq \tau(s) \land \tau(\sigma)$ and we let $\Perp$ be the set of all processes $t \star \pi$ such that $\tau(t) \land \tau(\pi) = 0$. 

\begin{theorem}
Let $\mathbb{B}$ be a complete Boolean algebra and $\delta$ a regular cardinal. $\mathbb{B}$ satisfies the $\delta$-cc if and only if $\mathcal{A}_\mathbb{B}$ satisfies the $\delta$-chain condition.
\end{theorem}

\begin{proof}
We begin by observing that any realizer is formed by using just $\lambda$-abstraction and the call-with-current-continuation and therefore for any realizer $\rlzfont{p} \in \mathcal{R}$, $\tau(\rlzfont{p}) = 1$. 

For the first direction, suppose that $\mathbb{B}$ satisfies the $\delta$-cc. We shall show that any realizer $\rlzfont{p}$ will witness that $\mathcal{A}_{\mathbb{B}}$ satisfies the $\delta$-chain condition. Let $\langle u_\beta \divline \beta \in \delta \rangle \subseteq \Lambda$ be a sequence of terms and suppose that for $\beta < \gamma$, $t \star u_\gamma \stackapp u_\beta \stackapp \pi \in \Perp$. By definition, this means that $0 = \tau(t) \land \tau(u_\beta) \land \tau(u_\gamma) \land \tau(\pi)$. Now, if $\tau(t) \land \tau(\pi) = 0$ then for any $\beta \in \delta$, $\tau(\rlzfont{p}) \land \tau(t) \land \tau(u_\beta) \land \tau(\pi) \in \Perp$ and therefore $\rlzfont{p} \star t \point u_\beta \point \pi \in \Perp$. So suppose that $\tau(t) \land \tau(\pi) > 0$. Then either $\tau(u_\beta) \land \tau(t) \land \tau(\pi) = 0$ for some $\beta \in \delta$ or $\tau(u_\gamma) \land \tau(u_\beta) = 0$ for every $\gamma \neq \beta$ in $\delta$. But, if the second case were to hold then $\{ \tau(u_\beta) \divline \beta \in \delta \}$ would form an antichain in $\mathbb{B}$, contradicting the assumption that $\mathbb{B}$ satisfies the $\delta$-cc. Thus we must be in the first case, in which case such a $u_\beta$ satisfies the conclusion of the claim.

For the reverse direction, suppose that $\mathcal{A}_\mathbb{B}$ satisfies the $\delta$-chain condition, as witnessed by $\rlzfont{p}$. Let $A \subseteq \mathbb{B}$ be an antichain of cardinality at least $\delta$ and enumerate the first $\delta$ many terms as $\langle p_\beta \divline \beta \in \delta \rangle$. For $p \in \mathbb{B}$, let $\omega_p$ be the stack bottom corresponding to $p$. Without loss of generality, we can assume that $0 \not\in A$. Then for every $a \neq b$ in $A$, we have that $\tau(\identity) \land \tau(\omega_a) \land \tau(\omega_b) \land \tau(\omega_1) = 0$ and therefore $\identity \star \omega_a \point \omega_b \point \omega_1 \in \Perp$. In particular, this gives us that if $\gamma < \beta$ then $\identity \star \omega_{p_\gamma} \stackapp \omega_{p_{\beta}} \stackapp \omega_1 \in \Perp$. Therefore, by the $\delta$-chain condition, there is some $\beta \in \delta$ such that $\rlzfont{p} \star \identity \point \omega_{p_\beta} \point \omega_1 \in \Perp$, which by definition means that $0 = \tau(\rlzfont{p}) \land \tau(\identity) \land \tau(\omega_{p_\beta}) \land \tau(\omega_1) = p_\beta$ which is a contradiction as we assumed that $a \neq 0$ for all $a \in A$.
\end{proof}

Note that since that realizability algebra is built from the Boolean algebra, it satisfies the following symmetric property: 
\Rcomment{$t \star s_0 \stackapp s_1 \stackapp \pi \in \Perp$ iff $t \star s_1 \stackapp s_0 \stackapp \pi \in \Perp$.} From this it follows that we can replace the $\delta$-chain condition with its uniform version.

\Rcomment{
\begin{proposition}
    Suppose that a realizability algebra $\mathcal{A}$ satisfies the following property:
    \[
    \forall t, s_0, s_1 \in \Lambda \; \forall \pi \in \Pi \; t \star s_0 \stackapp s_1 \stackapp \pi \in \Perp \Longleftrightarrow t \star s_1 \stackapp s_0 \stackapp \pi \in \Perp. 
    \]
    Then $\rlzfont{p}$ witnesses the $\delta$-chain condition for $\mathcal{A}$ iff it realizes the uniform $\delta$-chain condition for $\mathcal{A}$.
\end{proposition}

\begin{proof}
    We first see that the $\delta$-chain condition always implies the uniform version. To see this, suppose the $\delta$-chain condition holds, as witnessed by $\rlzfont{p}$, and let $A \subseteq \Lambda$ be a set of terms of cardinality at least $\delta$ such that for every $t \in \Lambda$ and $\pi \in \Pi$,
    \[
    \forall a, b \in A \; (a \neq b \rightarrow t \star a \stackapp b \stackapp \pi \in \Perp).
    \]
    Let $\langle u_\beta \divline \beta \in \delta \rangle$ be an enumeration of $\delta$ many of the terms. Then, whenever $\gamma < \beta$, $t \star u_\gamma \stackapp u_\beta \stackapp \pi \in \Perp$. Thus, by the $\delta$-chain condition, we can fix $\beta \in \delta$ for which $\rlzfont{p} \star t \stackapp u_\beta \stackapp \pi \in \Perp$. 

    Now suppose that $\rlzfont{p}$ witnesses the uniform $\delta$-chain condition and the additional property holds. Let $A = \langle u_\beta \divline \beta \in \delta \rangle$ be a sequence of terms such that for every $t \in \Lambda$ and $\pi \in \Pi$,
    \[
    \forall \gamma, \beta \in \delta (\gamma < \beta \rightarrow t \star u_\gamma \stackapp u_\beta \stackapp \pi \in \Perp).
    \]
    Then, if $u_\gamma$ and $u_\beta$ are in $A$ with $u_\gamma \neq u_\beta$, either $\gamma < \beta$ and $t \star u_\gamma \stackapp u_\beta \stackapp \pi \in \Perp$ or $\beta < \gamma$ and $t \star u_\beta \stackapp u_\gamma \stackapp \pi \in \Perp$. But, by the additional property, one of these processes is in the pole if and only if the other one is. Thus $t \star u_\gamma \stackapp u_\beta \stackapp \pi$ is always in $\Perp$. Hence, by the uniform $\delta$-chain condition, we can find some $u_\beta \in A$ such that $\rlzfont{p} \star t \stackapp u_\beta \stackapp \pi \in \Perp$. 
\end{proof}
}

\begin{corollary}
    Let $\mathbb{B}$ be a complete Boolean algebra and $\delta$ a regular cardinal. $\mathbb{B}$ satisfies the $\delta$-cc if and only if $\mathcal{A}_\mathbb{B}$ satisfies the \emph{uniform} $\delta$-chain condition.
\end{corollary}

\begin{definition}
    Let $\rlzin\text{-}\ff{Surj}(f, a, b)$ be an abbreviation for the \Rcomment{following} statement\Rcomment{:}
    \[
    \rlzin\text{-}\ff{Surj}(f, a, b) ~ \equiv ~ \forall y ( y \rlzin b \rightarrow \forall x (\op(x, y) \rlzin f \rightarrow x \notrlzin a) \rightarrow \perp)
    \]
\end{definition}

\Rcomment{Unpacking the above formula, $\rlzin\text{-}\ff{Surj}(f, a, b)$ is equivalent to $\forall y \rlzin b \exists x \rlzin a (\op(x, y) \rlzin f)$, which is the assertion that $f$ is an $\rlzin$-surjection from $a$ onto $b$.}

We begin with an observation:

\begin{proposition} \label{theorem:RealizingAlphaNeqBeta}
    Let $\rlzXNotInX \in \mathcal{R}$ be the fixed realizer such that $\rlzXNotInX \Vdash \forall x(x \not\in x)$ and let $\rlzXNotSimeqY = \lambda f \lambdaapp \lambda g \lambdaapp \app{g}{\rlzXNotInX}$. Then whenever $\alpha < \beta$, $\rlzXNotSimeqY \Vdash \fullname{\alpha} \not\simeq \fullname{\beta}$.
\end{proposition}

\begin{proof}
    Fix ordinals $\alpha < \beta$. Removing the abbreviation, we shall show that
    \[
    \lamAbstOne \lamAbstTwo \app{\lamtermTwo}{\rlzXNotInX} \Vdash \fullname{\alpha} \subseteq \fullname{\beta} \rightarrow (\fullname{\beta} \subseteq \fullname{\alpha} \rightarrow \perp).
    \]
    So, suppose that $t \Vdash \fullname{\alpha} \subseteq \fullname{\beta}$, $s \Vdash \fullname{\beta} \subseteq \fullname{\alpha}$ and take $\pi \in \Pi$. Then
    \[
    \rlzXNotInX \stackapp \pi \in \bigcup_{\gamma < \alpha} \{ r \stackapp \sigma \divline (\fullname{\gamma}, \sigma) \in \fullname{\beta}, r \Vdash \fullname{\gamma} \not\in \fullname{\alpha} \} = \falsity{ \fullname{\beta} \subseteq \fullname{\alpha}}.
    \]
    Thus, $\lamAbstOne \lamAbstTwo \app{\lamtermTwo}{\rlzXNotInX} \star t \stackapp s \stackapp \pi \succ s \star w \stackapp \pi \in \Perp$.
\end{proof}

We note again the anti-symmetric nature of this proof. Namely, if $\alpha < \beta$ then we need a different realizer to obtain $\rlzmodel \Vdash \fullname{\beta} \not\simeq \fullname{\alpha}$. This is precisely why we need to take an order for our chain of terms in the definition of $\delta$-chain conditions. 

\begin{theorem}
    Suppose that a realizability algebra satisfies the $\delta$-chain condition for some regular cardinal $\delta$, as witnessed by the term $\rlzfont{p}$. Then there exists a realizer $v$ such that,
    \[
    v \Vdash \forall f \forall a^{\fullname{\delta}} (\ff{ExtFun}(f, a) \rightarrow \rlzin\text{-}\ff{Surj}(f, a, \fullname{\delta}) \rightarrow \perp).
    \]
\end{theorem}

\begin{proof}
    Let $\rlzXNotSimeqY$ be the realizer from \Cref{theorem:RealizingAlphaNeqBeta}, let $\rlzNotNot = \lambda f \point \lambda g \point \app{g}{f}$ and let $\rlzXSimeqX \in \mathcal{R}$ be such that $\rlzXSimeqX \Vdash \forall x (x \simeq x)$. It is easy to see that, for any formula $\varphi$, $\rlzNotNot \Vdash \varphi \rightarrow (( \varphi \rightarrow \perp) \rightarrow \perp)$. We shall show that a realizer for the desired statement is
    \[
    v \coloneqq \lamAbst{t} \lamAbst{s} 
    \Bigglapp{
    \biggrapp{\cc}{
    \lamAbst{k} 
    \Bigtwoapp{
    \bigrapp{\rlzNotNot}{\lamAbstOne \bigrapp{k}{\biglapp{\rapp{\rlzfont{p}}{\app{t}{\rlzXSimeqX}}}{\lamtermOne}}}
    }
    {\app{s}{\identity}}
    }
    }
    {\rlzXNotSimeqY}.
    \]

    Due to the complicated nature of this realizer, we briefly explain where the construction comes from. The idea it to suppose that $v$ did not realize our statement. This means that we will be able to find a set $f$, ordinal $\gamma \in \delta$, stack $\pi$ and terms $t, s$ such that $t$ realizes $f$ is an extensional function with domain $\fullname{\gamma}$ and $s$ realizes $f$ is an $\rlzin$-Surjection but $v \star t \stackapp s \stackapp \pi \not\in \Perp$. Deconstructing $v$ will give a new term $v'$ such that for any $\beta \in \delta$
    \[
    v' \not\Vdash \forall x (\op(x, \fullname{\beta}) \rlzin f \rightarrow x \notrlzin \fullname{\gamma}).
    \]
    Therefore, for each $\beta \in \delta$ we can choose some $\alpha_\beta \in \gamma$, stack $\pi_\beta \in \Pi$ and term $u_\beta$ such that $u_\beta \Vdash \op(\fullname{\alpha_\beta}, \fullname{\beta}) \rlzin f$ but $v' \star u_\beta \stackapp \pi_\beta \not\in \Perp$. Simplifying $v'$ again will then give us that for each $\beta \in \delta$, $\rlzfont{p} \star \app{t}{\rlzXSimeqX} \stackapp u_\beta \stackapp \rlzXNotSimeqY \stackapp \pi \not\in \Perp$.

    However, as $\delta$ is regular and $\gamma \in \delta$, we can find some $\alpha \in \gamma$ for which $B \coloneqq \{ \beta \in \delta | \alpha_\beta = \alpha \}$ has cardinality $\delta$. Then, since the $u_\beta$ realize incompatible information for $f(\fullname{\alpha})$ whenever $\beta < \beta'$, we will have $\app{t}{\rlzXSimeqX} \star u_\beta \stackapp u_{\beta'} \stackapp \rlzXNotSimeqY \stackapp \pi \in \Perp$ for all $\beta < \beta'$ in $B$. But this will contradict the $\delta$-chain condition.

    \bigskip
    
   \noindent  We now give the details of the argument.
    
    Suppose for a contradiction that $v$ did not realize our statement. Then we could fix some $f \in \rlzstr$, $\gamma \in \delta$, $t \Vdash \ff{ExtFun}(f, \fullname{\gamma})$, $s \Vdash \rlzin\text{-}\ff{Surj}(f, \fullname{\gamma}, \fullname{\delta})$ and $\pi \in \Pi$ such that $v \star t \stackapp s \stackapp \pi \not\in \Perp$. Now this means that
    \begin{align*}
        v \star t \point s \point \pi & \succ \Bigglapp{
    \biggrapp{\cc}{
    \lamAbst{k} 
    \Bigtwoapp{
    \bigrapp{\rlzNotNot}{\lamAbstOne \bigrapp{k}{\biglapp{\rapp{\rlzfont{p}}{\app{t}{\rlzXSimeqX}}}{\lamtermOne}}}
    }
    {\app{s}{\identity}}
    }
    }
    {\rlzXNotSimeqY}
    \star \pi \\
    & \succ \cc \star 
    \Big(\lamAbst{k} 
    \Bigtwoapp{
    \bigrapp{\rlzNotNot}{\lamAbstOne \bigrapp{k}{\biglapp{\rapp{\rlzfont{p}}{\app{t}{\rlzXSimeqX}}}{\lamtermOne}}}
    }
    {\app{s}{\identity}}\Big)
    \stackapp \rlzXNotSimeqY \stackapp \pi \\
    & \succ 
    \lamAbst{k} 
    \Big(\bigtwoapp{
    \bigrapp{\rlzNotNot}{\lamAbstOne \bigrapp{k}{\biglapp{\rapp{\rlzfont{p}}{\app{t}{\rlzXSimeqX}}}{\lamtermOne}}}
    }
    {\app{s}{\identity}}\Big)
    \star \saverlz{\rlzXNotSimeqY \stackapp \pi} \stackapp \rlzXNotSimeqY \stackapp \pi \\
    & \succ
    \bigrapp{\rlzNotNot}{\lamAbstOne \bigrapp{\saverlz{\rlzXNotSimeqY \stackapp \pi}}{\biglapp{\rapp{\rlzfont{p}}{\app{t}{\rlzXSimeqX}}}{\lamtermOne}}}
    \star \app{s}{\identity} \stackapp \rlzXNotSimeqY \stackapp \pi.
    \end{align*}
    So, in particular, we must have that $\bigrapp{\rlzNotNot}{\lamAbstOne \bigrapp{\saverlz{\rlzXNotSimeqY \stackapp \pi}}{\biglapp{\rapp{\rlzfont{p}}{\app{t}{\rlzXSimeqX}}}{\lamtermOne}}}
    \star \app{s}{\identity} \stackapp \rlzXNotSimeqY \stackapp \pi \not\in \Perp$.
    
    Next, we observe that for any $\beta \in \delta$, $\identity \Vdash \fullname{\beta} \rlzin \fullname{\delta}$. Therefore, since $s$ realizes that $f$ is a surjection from $\fullname{\gamma}$ onto $\fullname{\delta}$, for every $\beta \in \delta$ we have
    \[
    \app{s}{\identity} \Vdash \forall x (\op(x, \fullname{\beta}) \rlzin f \rightarrow x \notrlzin \fullname{\gamma}) \rightarrow \perp.
    \]
    From which we can deduce that $\app{s}{\identity} \stackapp \rlzXNotSimeqY \stackapp \pi \in \falsity{(\forall x (\op(x, \fullname{\beta}) \rlzin f \rightarrow x \notrlzin \fullname{\gamma}) \rightarrow \perp) \rightarrow \perp}$. Thus, for every $\beta \in \delta$ we must have that 
    \[
    \bigrapp{\rlzNotNot}{\lamAbstOne \bigrapp{\saverlz{\rlzXNotSimeqY \stackapp \pi}}{\biglapp{\rapp{\rlzfont{p}}{\app{t}{\rlzXSimeqX}}}{\lamtermOne}}}
    \not\Vdash (\forall x (\op(x, \fullname{\beta}) \rlzin f \rightarrow x \notrlzin \fullname{\gamma}) \rightarrow \perp) \rightarrow \perp.
    \]
    Since $\rlzNotNot$ was constructed to be a realizer for double negations this gives us that
    \[
    \lamAbstOne \bigrapp{\saverlz{\rlzXNotSimeqY \stackapp \pi}}{\biglapp{\rapp{\rlzfont{p}}{\app{t}{\rlzXSimeqX}}}{\lamtermOne}}
    \not\Vdash \forall x (\op(x, \fullname{\beta}) \rlzin f \rightarrow x \notrlzin \fullname{\gamma}).
    \]
    Deconstructing what that means, for every $\beta \in \delta$ we can choose an $\alpha_\beta \in \gamma$, $\pi_\beta \in \Pi$ and \linebreak[4] $u_\beta \Vdash \op(\fullname{\alpha_\beta}, \fullname{\beta}) \rlzin f$ such that 
    $\lamAbstOne \bigrapp{\saverlz{\rlzXNotSimeqY \stackapp \pi}}{\biglapp{\rapp{\rlzfont{p}}{\app{t}{\rlzXSimeqX}}}{\lamtermOne}} 
    \star u_\beta \stackapp \pi_\beta \not\in \Perp$. 
    Again simplifying the term, we have that for every $\beta \in \delta$
    \begin{align*}
        \lamAbstOne \bigrapp{\saverlz{\rlzXNotSimeqY \stackapp \pi}}{\biglapp{\rapp{\rlzfont{p}}{\app{t}{\rlzXSimeqX}}}{\lamtermOne}} \star u_\beta \stackapp \pi_\beta & \succ
        \saverlz{\rlzXNotSimeqY \stackapp \pi} \star \biglapp{\rapp{\rlzfont{p}}{\app{t}{\rlzXSimeqX}}}{u_\beta} \stackapp \pi_\beta \\
        & \succ
        \biglapp{\rapp{\rlzfont{p}}{\app{t}{\rlzXSimeqX}}}{u_\beta} \star \rlzXNotSimeqY \stackapp \pi \\
        & \succ 
        \rlzfont{p} \star \app{t}{\rlzXSimeqX} \stackapp u_\beta \stackapp \rlzXNotSimeqY \stackapp \pi.
    \end{align*}
    So, for all $\beta \in \delta$, $\rlzfont{p} \star \app{t}{\rlzXSimeqX} \stackapp u_\beta \stackapp \rlzXNotSimeqY \stackapp \pi \not\in \Perp$.

    However, the function $\beta \mapsto \alpha_\beta$ maps $\delta$ to $\gamma$. Therefore, by the regularity of $\delta$ there must exist an $\alpha \in \gamma$ for which $B \coloneqq \{ \beta \in \delta \divline \alpha_\beta = \alpha \}$ has cardinality $\delta$. Then, since $u_\beta \Vdash \op(\fullname{\alpha}, \fullname{\beta}) \rlzin f$ and $t$ realizes that $f$ is an extensional function, if $\beta \in \beta'$ then 
    \begin{multline*}
    \rlzXSimeqX \stackapp u_\beta \stackapp u_{\beta'} \stackapp \rlzXNotSimeqY \stackapp \pi \\
    \in \falsity{ \fullname{\alpha} \simeq \fullname{\alpha} \rightarrow \op(\fullname{\alpha}, \fullname{\beta}) \rlzin f \rightarrow \op(\fullname{\alpha}, \fullname{\beta'}) \rlzin f \rightarrow \fullname{\beta} \not\simeq \fullname{\beta'} \rightarrow \perp}. 
    \end{multline*}
    So, $\app{t}{\rlzXSimeqX} \star u_\beta \stackapp u_{\beta'} \stackapp \rlzXNotSimeqY \stackapp \pi \in \Perp$ whenever $\beta < \beta'$. Finally, by the $\delta$-chain condition, this means that we can fix some $ \beta \in B$ such that $\rlzfont{p} \star \app{t}{\rlzXSimeqX} \stackapp u_\beta \stackapp \rlzXNotSimeqY \stackapp \pi \in \Perp$, yielding the desired contradiction. 
\end{proof}

\begin{corollary}
    Suppose that a realizability algebra satisfies the $\delta$-chain condition for some regular cardinal $\delta$, as witnessed by the term $\rlzfont{p}$. Then $\fullname{\delta}$ is an $\rlzin$-cardinal in $\rlzmodel$ and hence also an $\in$-cardinal.
\end{corollary}

Recall that $\rlzin\text{-}\ff{Fun}(f, a)$ abbreviates the statement that $f$ is an $\rlzin$-function with domain $a$. That is
\[
\rlzin\text{-}\ff{Fun}(f, a)  ~ \equiv ~ \forall x \forall y \forall y' (\op(x, y) \rlzin f \rightarrow \op(x, y') \rlzin f \rightarrow y \neq y' \rightarrow \perp).
\]
If our realizability algebra only satisfies the uniform $\delta$-chain condition then we can only prove that for every $a \rlzin \fullname{\delta}$, if $\rlzin\text{-}\ff{Fun}(f, a)$ then $f$ is not an $\rlzin$-Surjection. This is because in the above proof, if $t \Vdash \ff{ExtFun}(f, a)$ then for $\beta \neq \beta'$ we will have
\[
t \rlzXSimeqX u_\beta u_{\beta'} \Vdash \fullname{\beta} \not\simeq \fullname{\beta'} \rightarrow \perp.
\]
However, we cannot deduce $t \rlzXSimeqX \star u_\beta \stackapp u_{\beta'} \stackapp \pi \in \Perp$ from this because we need different terms depending on whether $\beta < \beta'$ or $\beta' < \beta$. This problem is resolved when $t \Vdash \rlzin\text{-}\ff{Fun}(f, \fullname{\gamma})$ due to the symmetric nature of non-extensionality equality.

We further remark that if it were to be the case that for all $\beta, \beta' \in B$, $u_\beta = u_{\beta'}$, then in order to apply the uniform $\delta$-chain condition we formally need to define distinct terms which realize the same things. This can be easily done. For example, let $\langle \nu_\alpha \divline \alpha \in |\Lambda| \rangle$ be an enumeration of $\Lambda$. Then we can define $\bar{u}_\beta = \rapp{\lamAbstTwo \lamAbstThree \lamtermTwo}{\app{u_\beta}{\nu_\beta}}$. It is then easy to see that 
\[
u_\beta \Vdash \varphi \quad \Rightarrow \quad \bar{u}_\beta \Vdash \varphi. 
\]

\begin{theorem}
    Suppose that a realizability algebra satisfies the uniform $\delta$-chain condition for some regular cardinal $\delta$, as witnessed by the term $\rlzfont{p}$. Then there exists a realizer $v$ such that,
    \[
    v \Vdash \forall f \forall a^{\fullname{\delta}} (\rlzin\text{-}\ff{Fun}(f, a) \rightarrow \rlzin\text{-}\ff{Surj}(f, a, \fullname{\delta}) \rightarrow \perp).
    \]
\end{theorem}

\bibliography{main}
\bibliographystyle{alpha}

\end{document}